\documentclass[reqno]{amsart}
\usepackage{hyperref}
\usepackage{mathrsfs}

\AtBeginDocument{}

\begin{document}
\title[\hfil Regularity of pullback attractors and equilibria]
{Regularity of pullback attractors and equilibria for a stochastic non-autonomous reaction-diffusion equations perturbed by a multiplicative noise}

\author{Wenqiang Zhao, Shuzhi Song}  

\address{Wenqiang Zhao, Shuzhi Song \newline
School Of Mathematics and Statistics, Chongqing Technology and Business
University,  Chongqing 400067, China}
\email{gshzhao@sina.com}

\keywords{Random dynamical systems;
non-autonomous reaction-diffusion \hfill\break\indent equation;
 upper semi-continuity; $\mathcal{D}$-pullback attractor; multiplicative noise}

\begin{abstract}
In this paper, a standard about the existence and upper semi-continuity of pullback attractors  in the non-initial space is established for
some classes of non-autonomous SPDE. This pullback attractor, which is the omega-limit set of the absorbing set constructed in the initial space, is completely determined
by the asymptotic compactness of solutions in both the initial and non-initial spaces. As applications,
the existences and upper semi-continuity  of pullback attractors in $H^1(\mathbb{R}^N)$ are proved
for stochastic non-autonomous reaction-diffusion equation  driven by a multiplicative noise.
Finally we show that under some additional conditions the cocycle  admits a unique equilibrium.
\end{abstract}

\maketitle
\numberwithin{equation}{section}
\newtheorem{theorem}{Theorem}[section]
\newtheorem{lemma}[theorem]{Lemma}
\newtheorem{remark}[theorem]{Remark}
\newtheorem{definition}[theorem]{Definition}
\allowdisplaybreaks

\section{Introduction}

In this paper, we consider the dynamics of solutions to
the following reaction-diffusion equation on $\mathbb{R}^N$ driven by a random noise as well as  a deterministic non-autonomous forcing:
\begin{align} \label{eq1}
du+(\lambda u-\Delta u)dt=f(x,u)dt+g(t,x)dt+\varepsilon u\circ d\omega(t),
\end{align}
with initial  condition
\begin{align} \label{eq2}
u(\tau,x)=u_0(x), \ \  \ x\in\mathbb{R}^N,
\end{align}
where the initial $u_0\in L^2(\mathbb{R}^N)$, $\lambda$ is a positive constant, $\varepsilon$ is the intensity of noise, the unknown $u=u(x,t)$ is a
real valued function of $x\in \mathbb{R}^N$ and $t> \tau$, $\omega(t)$ is a mutually independent
two-sided real-valued Wiener process defined on a
 canonical Wiener probability space $(\Omega,
\mathcal {F}, {P})$.

The notion of random attractor, introduced in \cite{Rand1,Rand2,Rand3,Rand4}, is
an important tool to study the qualitative property of stochastic partial differential equations(SPDE). We can find  a large number of literature to investigate
the existences of random attractors in \emph{an initial space} (the initial data located space) for  some concrete stochastic partial differential equations, see \cite{Bates, And,Wang1,Wang4, Wangz, Zhang, Zhao} and the references
therein.  In particular, \cite{Wang3,Wang1,Wang5} discussed the upper semi-continuity of a family of random attractors in the initial space.

As we know, the solutions of SPDE may possess some regularity, for example, higher-order integrability or higher-order differentiability. In these cases, the the solutions  may escape (or leave) the initial space and enter into another space, which we call \emph{a non-initial space}.
So, the existence and upper semi-continuity of random attractors in a non-initial space, usually a higher-regularity space, such as $L^p$($p>2$) or
$H^1$, are necessary for us  to understand the dynamics of solutions of SPDE.

In terms of  this consideration, some literature attacked this problem recently. In the case of bounded domain, \cite{Cun, Lijia,Liyangrong2,Zhao1,Zhao2,Yin}
discussed the existence of random attractor in the non-initial spaces $L^p$ and $H^1_0$ space, respectively.
When the state space is unbounded,  Zhao and Li \cite{Zhao3} proved the existence of random attractors for reaction-diffusion equations
in $L^p(\mathbb{R}^N)$,  and for the same equation, Li and et al \cite{Liyangrong1} obtained the upper semi-continuity of random attractor
in $L^p(\mathbb{R}^N)$.
Most recently, Zhao \cite{Zhao0, Zhao4} proved the existence of random attractors for semi-linear degenerate parabolic equations in $L^{2p-2}(D)\cap H^1(D)$, where $D$ is a unbounded domain. Bao \cite{Bao} proved the existence of random attractors for non-autonomous Fitzhugh-Nagumo system  in $H^1(\mathbb{R}^N)\times L^2(\mathbb{R}^N)$. In that paper, the key point is closely related to Lemma 5.1 there, of which the detailed proof is omitted.

It is pointed out that most recently, Li and et al \cite{Liyangrong3} established the theory of bi-spatial random attractors by using the notion of uniform omega-limit compactness, by which SPDE with autonomous forcing can be solved, see also \cite{Liyangrong1}. However, to the best of our knowledge, there are
 no literature to discuss the existence and  upper semi-continuity of pullback attractors in a non-initial space for SPDE with a non-autonomous forcing term, except  the literature \cite{Bao}.

  In this paper, we study the existence and upper semi-continuity of pullback attractors in the non-initial space $H^1(\mathbb{R}^N)$ for problem (\ref{eq1})-(\ref{eq2}) with a non-autonomous forcing. The nonlinearity $f$ and the deterministic non-autonomous function $g$ satisfy almost the same conditions as \cite{Wang1},
in which the author obtained the existence and upper-continuity of pullback attractors in the initial space $L^2(\mathbb{R}^N)$.  Here, we strengthen this result and show that
the obtained pullback attractors are also compact and attracting in $H^1(\mathbb{R}^N)$ norm.  Furthermore, we find that the upper continuity of the obtained pullback attractors happen
in $H^1(\mathbb{R}^N)$.
The existence of pullback attractor in an initial space  for a non-autonomous stochastic partial differential equation is
established in \cite{Wang4}, where the measurability of pullback attractors is proved. The applications we may see \cite{And,Wang1,Wang3,Wang4}.
For the reference on the theory regarding  upper semi-continuity of pullback attractors, we may refer to \cite{Wang1,Wang3,Wang5} for the stochastic cases and to \cite{James,Hale} for the deterministic cases.

In order to solve our problem, we establish a sufficient criteria for
the existence and upper semi-continuity of pullback attractors in a non-initial space. It is showed that a family of pullback attractors obtained in an initial space are compact, attracting and upper semi-continuous in a non-initial space
 if  some compactness conditions of the cocycles are satisfied, see Theorem 2.6-2.8 in section 2. This implies that
 the continuity (or quasi-continuity \cite{Liyangrong2}, norm-weak continuity \cite{Zhong}) and absorption in the non-initial space are not necessary ones.
  This result is a  meaningful and convenient tool for us to consider the  existence and
 upper semi-continuity of pullback attractors in some associated non-initial spaces for SPDE with a non-autonomous forcing term.

 Consider that the stochastic equation (\ref{eq1}) is defined on unbounded domains, the asymptotic compactness of solution in $H^1(\mathbb{R}^N)$ can not be
derived by the traditional technique. The reasons are as follows. On the one hand, the equation (\ref{eq1}) is stochastic and the Wiener process $\omega$ is only continuous but not
differentiable in $t$. This leads to some difficulties for us to estimate the norm of derivative  $u_t$ by the trick  employed in deterministic case, see \cite{Zhong,ZhangY}. Then
the asymptotic
compactness in $H^1(\mathbb{R}^N)$ can not be proved by estimate of the difference of $\nabla u$  as in \cite{ZhangY}.

On the other hand,
the estimate of $\Delta u$ is not available for our problem (up to now, actually we do not know how to estimate the norm $\Delta u$ of problem (\ref{eq1})-(\ref{eq2}), although this can be achieved in deterministic case by estimate $u_t$, see \cite{Zhong}). So we can not use the Sobolev compact embeddings of $H^2\hookrightarrow H^1$ on bounded domains.
 Here we surmount these obstacles by checking the uniform smallness of solutions outside a large ball in $H^1(\mathbb{R}^N)$-norms as in \cite{Wang2,Wang3,Bao}. In  bounded domains, we prove the asymptotic compactness of solutions by a space splitting technique as in \cite{Zhao0,Zhao1,Zhao2, Cun}, and combination the estimate of the truncation of solutions in $L^{2p-2}$-norm over an integral interval.

 Finally, we investigate how the parameters in problem (1.1) affect the pullback attractor. We show that if the parameters satisfy some conditions, then the cocycle admits a unique
 equilibrium in $L^2(\mathbb{R}^N)$. Furthermore, this  equilibrium is also in both $H^1(\mathbb{R}^N)$ and $L^p(\mathbb{R}^N)$, $p>2$.

In the next section, we recall some notions  and prove a sufficient standard for the existence and upper semi-continuity of pullback attractors of
non-autonomous system in a non-initial space. In section 3, we give the assumptions on $g$ and $f$, and define a continuous cocycle for problem (\ref{eq1})-(\ref{eq2}).
In section 4 and 5,  we prove the existence and upper semi-continuity in $H^1(\mathbb{R}^N)$. Finally,  in section 6,
we prove the existence of equilibria for the cocyle derived from problem (1.1).

\section{Preliminaries and abstract results}

Let $(X, \|.\|_X)$ and $(Y, \|.\|_Y)$  be two complete separable Banach spaces  with Borel
sigma-algebras $\mathcal {B}(X)$ and $\mathcal {B}(Y)$, respectively. $X\cap Y\neq \emptyset$. For convenience, we call $X$ the initial space ( which contains all initial data of a SPDE) and $Y$ the associated non-initial space (usually the  regular solutions (of a SPDE) located space).

In this section,  we give a sufficient standard for the existence and upper semi-continuity  of pullback attractors
in the non-initial  space $Y$ for
the random dynamical system (RDS) over two parametric spaces.  The readers may refer to \cite{Zhao0,Zhao3,Zhao4,Lijia,Liyangrong1,Liyangrong2,Liyangrong3,Yin} for the existence and semi-continuity of random attractors in the non-initial  space $Y$ for a RDS over one parametric space.

We also mention that regarding the existence of random attractors  in the initial space $X$ for the RDS over one parametric space, the good references
 are \cite{Arn,Bates,Rand1,Rand2,Rand3,Rand4}. However, here we recall from  \cite{Wang4} some basic notions regarding RDS over two parametric spaces,  one of which is the
 real numbers space and another of which is  the measurable probability space with a measure preserving transformation.

\subsection{Preliminaries}

The basic notion in RDS is a metric (or measurable)
dynamical system (MDS)
$\vartheta\equiv(\Omega,\mathcal{F},{P},\{\vartheta_t\}_{t\in
\mathbb{R}})$, which is  a probability space
$(\Omega,\mathcal{F},{P})$ with a group $\vartheta_t,t\in
\mathbb{R},$ of measure preserving transformations of
$(\Omega,\mathcal{F},{P})$.

An MDS $\vartheta$ is said to be
ergodic under ${P}$ if for any $\vartheta$-invariant set $F\in
\mathcal {F}$, we have either ${P}(F)=0$ or ${P}(F)=1$,
where the $\vartheta$-invariant set is in the sense that  ${P}(\vartheta_tF)=(F)$ for
$F\in \mathcal {F}$ and all $t\in \mathbb{R}$.\\

\textbf{Definition 2.1}.
\emph{Let $(\Omega,\mathcal {F},{P},\{\vartheta_t\}_{t\in\mathbb{R}})$ be a measurable dynamical system. A family of measurable mappings
 $\varphi: \mathbb{R}^+\times\mathbb{R}\times \Omega\times X\mapsto X$ is called  a cocycle on $X$ over $\mathbb{R}$ and $(\Omega,\mathcal {F},{P},\{\vartheta_t\}_{t\in\mathbb{R}})$
if for all $\tau\in \mathbb{R}, \omega\in\Omega$ and $t,s\in \mathbb{R}^+$, the following conditions are satisfied:
\begin{align}
&\varphi(0, \tau,\omega,.)\ \mbox{is the identity on X},\notag\\
&\varphi(t+s,\tau, \omega,.)=\varphi(t,\tau+s, \vartheta_s\omega, .)\circ \varphi(s,\tau, \omega, .).\notag
 \end{align}
In addition, if $\varphi(t,\tau, \omega,.): X\rightarrow X$ is continuous for all $t\in
\mathbb{R}^+, \tau\in \mathbb{R}, \omega\in\Omega$, then $\varphi$ is called a continuous cocycle on $X$ over $\mathbb{R}$ and $(\Omega,\mathcal {F},{P},\{\vartheta_t\}_{t\in\mathbb{R}})$.}\\

\textbf{Definition 2.2.}\emph{ Let $2^X$ be the collection of all subsets of $X$. A set-valued mapping $K: \mathbb{R}\times \Omega\mapsto 2^X$ is called measurable in $X$ with respect to $\mathcal{F}$ in $\Omega$ if  the mapping $\omega\in\Omega\mapsto \mbox{dist}_X(x,K(\tau,\omega))$ is ($\mathcal{F},\mathcal{B}(\mathbb{R}))$-measurable for every fixed $x\in X$
and $\tau\in\mathbb{R}$, where $\mbox{dist}_X$ is the Haustorff semi-metric in $X$. In this case, we also say the family $\{K(\tau,\omega);\tau\in\mathbb{R},\omega\in\Omega\}$ is measurable in $X$ with respect to $\mathcal{F}$ in $\Omega$. Furthermore if the value
$K(\tau,\omega)$ is a closed nonempty subset of $X$ for all $\tau\in\mathbb{R}$ and $\omega\in \Omega$, then $\{K(\tau,\omega);\tau\in\mathbb{R},\omega\in\Omega\}$ is called a closed measurable set of $X$ with respect to $\mathcal{F}$ in $\Omega$.}\\

Hereafter, we always assume that $\varphi$ is a continuous cocycle on $X$ over $\mathbb{R}$ and $(\Omega,\mathcal {F},{P},\{\vartheta_t\}_{t\in\mathbb{R}})$ satisfying

\emph{(H1)£º\ \  For every fixed $t\in\mathbb{R}^+, \tau\in\mathbb{R}$ and $\omega\in\Omega$, $\varphi(t,\tau,\omega,.): X\mapsto X\cap Y$; }\\

\emph{(H2)\ \  If $\{x_n\}_n\subset X\cap Y$ such that $x_n\rightarrow x$ in
$X$ and $x_n\rightarrow y$ in $Y$, respectively, then $x=y$.}\\

Let $\mathcal{D}$ be  a  collection  of some families of nonempty subsets of $X$ parametrized by $\tau\in\mathbb{R}$ and $\omega\in \Omega$ such that
$$
\mathcal{D}=\{B=\{B(\tau,\omega)\in 2^X;B(\tau,\omega)\neq\emptyset,\tau\in\mathbb{R},\omega\in \Omega \};f_B\ \mbox{satisfies some conditions}\}.
$$
In particular, for two elements $B_1,B_2\in\mathcal{D}$, we say that $B_1=B_2$ if and only if $B_1(\tau,\omega)=B_2(\tau,\omega)$ for all $\tau\in\mathbb{R}$ and $\omega\in \Omega$.\\

\textbf{Definition 2.3.} \emph{Let  $\mathcal{D}$ be  a  collection  of some families of nonempty subsets of $X$ and $K=\{K(\tau,\omega);\tau\in\mathbb{R},\omega\in\Omega\}\in \mathcal{D}$. Then $K$ is
  called a $\mathcal{D}$-pullback absorbing set for
a cocycle  $\varphi$ in $X$ if all $\tau\in \mathbb{R}, \omega\in\Omega$ and for every $B\in \mathcal{D}$
 there exists a absorbing time $T=T(\tau,\omega,B)>0$ such that
$$\varphi (t,\tau-t, \vartheta_{-t}\omega, B(\tau-t,\vartheta_{-t}\omega))\subseteq K(\tau, \omega)\ \ \ \ \ \ \mbox{for all}\ t\geq
T.$$
If in addition, $K$ is measurable in $X$ with respect to the $P$-completion of $\mathcal{F}$ in $\Omega$, then $K$ is said to a  measurable pullback absorbing set
for $\varphi$.}\\

 \textbf{Definition 2.4.} \emph{Let  $\mathcal{D}$ be  a  collection  of some families of nonempty subsets of $X$. A cocycle  $\varphi$ is said to be
$\mathcal{D}$-pullback asymptotically compact in $X(\mbox{resp. in}\ Y)$ if for all $\tau\in \mathbb{R}, \omega\in\Omega$
$$\{\varphi(t_n, \tau-t_n, \vartheta_{-t_n}\omega, x_n)\}\ \mbox{has a
convergent subsequence in}\ X(\mbox{resp. in}\ Y)$$
whenever $t_n\rightarrow \infty$ and
$x_n\in B(\tau-t_n,\vartheta_{-t_n}\omega)$ with $B=\{B(\tau, \omega);\tau\in \mathbb{R},\omega\in \Omega\} \in \mathcal{D}$.}\\

\textbf{Definition 2.5.} \emph{ Let  $\mathcal{D}$ be  a  collection  of some families of nonempty subsets of $X$ and $\mathcal{A}=\{\mathcal{A}(\tau,\omega);\tau\in\mathbb{R},\omega\in\Omega\}\in \mathcal{D}$. $\mathcal{A}$ is called  a $\mathcal{D}
$-pullback  attractor for a cocycle $\varphi$ in $X(\mbox{resp. in}\ Y)$ over $\mathbb{R}$ and $(\Omega,\mathcal{F},P,\{\vartheta_t\}_{t\in\mathbb{R}})$ if
}\\
\emph{(i)  $\mathcal {A}$ is measurable in $X$ with respect to the $P$-completion of $\mathcal{F}$, and $\mathcal {A}(\tau,\omega)$ is compact in $X\ (\mbox{resp. in}\ Y)$ for all $\tau\in \mathbb{R}, \omega\in\Omega$;
}\\
\emph{(ii) $\mathcal {A}$ is invariant, that is, for all $\tau\in \mathbb{R}, \omega\in\Omega$,
 $$
 \varphi(t,\tau,\omega,\mathcal{A}(\tau,\omega))=\mathcal{A}(\tau+t,\vartheta_t\omega), \forall\ t\geq0;
 $$}
\emph{(iii) $\mathcal {A}$  attracts every element $B=\{B(\tau, \omega);\tau\in \mathbb{R},\omega\in \Omega\} \in \mathcal{D}$ in $X(\mbox{resp. in}\ Y)$, that is, for all $\tau\in \mathbb{R}, \omega\in\Omega$,
$$\mathop {\lim }\limits_{t\rightarrow +\infty}\mbox{dist}_X(\varphi(t,\tau-t,\vartheta_{-t}\omega, B(\tau-t, \vartheta_{-t}\omega)),\mathcal{A}(\tau,\omega))= 0$$
$$(\mbox{resp.}\mathop {\lim }\limits_{t\rightarrow +\infty}\mbox{dist}_Y(\varphi(t,\tau-t,\vartheta_{-t}\omega, B(\tau-t, \vartheta_{-t}\omega)),\mathcal{A}(\tau,\omega))= 0).$$
}

\subsection{Existence of random attractors  in a non-initial space $Y$}

 This subsection is concerned
with the existence of $\mathcal{D}$-pullback attractor of the cocycle $\varphi$ in the non-initial space $Y$. The continuity of $\varphi$ in $Y$ is not clear, and the inclusion relation of
$X$ and $Y$ is also unknown except that $\emph{(H1)}$ hand $\emph{(H2)}$ hold.\\

 \textbf{ Theorem  2.6. }  \emph{Let $\mathcal{D}$ be a collection of some families of nonempty subsets of $X$ which is inclusion closed. Let $\varphi$ be a continuous cocycle on $X$ over $\mathbb{R}$ and $(\Omega,\mathcal{F},P,\{\vartheta_t\}_{t\in\mathbb{R}})$. Assume that
}\\
\emph{(i) $\varphi$ has  a closed and measurable (w.r.t. the the $P$-completion of $\mathcal{F}$) $\mathcal{D}$-pullback absorbing set
$K=\{K(\tau,\omega);\tau\in\mathbb{R},\omega\in\Omega\}\in \mathcal{D}$ in $X$;
}\\
\emph{(ii) $\varphi$ is  $\mathcal{D}$-pullback asymptotically compact in $X$.}
\emph{Then the cocycle  $\varphi$ has a unique $\mathcal{D}$-pullback  attractor $\mathcal{A}_X=\{\mathcal{A}_X(\tau,\omega); \tau\in\mathbb{R}, \omega\in \Omega\}\in\mathcal{D}$ in $X$,
given by
\begin{align}\label{ff01}
\mathcal{A}_X(\tau,\omega)=\bigcap_{s>0}\overline{\bigcup_{t\geq s} \varphi(t,\tau-t,\vartheta_{-t}\omega, K(\tau-t,\vartheta_{-t}\omega))}^X, \ \
\tau\in \mathbb{R},\omega\in \Omega,
 \end{align} where the closure is taken in $X$.}

\emph{ If further
(H1)-(H2) hold  and (iii) $\varphi$ is  $\mathcal{D}$-pullback asymptotically compact in $Y$,}
\emph{ then the cocycle  $\varphi$ has a unique $\mathcal{D}$-pullback  attractor $\mathcal{A}_Y=\{\mathcal{A}_Y(\tau,\omega); \tau\in\mathbb{R}, \omega\in \Omega\}$ in $Y$,
given by
\begin{align}\label{ff02}
\mathcal{A}_Y(\tau,\omega)=\bigcap_{s>0}\overline{\bigcup_{t\geq s} \varphi(t,\tau-t,\vartheta_{-t}\omega, K(\tau-t,\vartheta_{-t}\omega))}^Y, \ \
\tau\in \mathbb{R},\omega\in \Omega.
 \end{align}}

\emph{In addition, we have $\mathcal{A}_Y=\mathcal{A}_X\subset X\cap Y$ in the sense of set inclusion, i.e.,
 for every $\tau\in \mathbb{R},\omega\in \Omega$, $\mathcal{A}_Y(\tau,\omega)=\mathcal{A}_X(\tau,\omega)$.}\\

\emph{ Proof} \ \  The first result is well known and so we are interested in the second result.
Indeed, (\ref{ff02}) makes sense by (\emph{H1}) and $\mathcal{A}_Y\neq \emptyset$
by the asymptotic compactness of the cocycle $\varphi$ in $Y$. In the following, we  show that  $\mathcal{A}_Y$
satisfies  Definition 2.5.

 \emph{Step 1.} We claim that the set $\mathcal{A}_Y$ is measurable  in $X$ (w.r.t the $P$-completion of $\mathcal{F}$ in $\Omega$) and $\mathcal{A}_Y\in\mathcal{D}$ is invariant  by proving that $\mathcal{A}_Y=\mathcal{A}_X$ since $\mathcal{A}_X$ is measurable (with respect to the $P$-completion of $\mathcal{F}$ in $\Omega$)  and $\mathcal{A}_X\in\mathcal{D}$ is invariant (the measurability of $\mathcal{A}_X$ is proved by Theorem 2.14 in \cite{Wang3}).

For every fixed $\tau\in\mathbb{R}$ and $\omega\in\Omega$, taking $x\in \mathcal{A}_X(\tau,\omega)$, by (\ref{ff01}), there exist two sequences $t_n\rightarrow+\infty$
and $x_n\in K(\tau-t_n,\vartheta_{-t_n}\omega)$ such that
\begin{align} \label{ff3}
 \varphi(t_n, \tau-t_n, \vartheta_{-t_n}\omega, x_n)\xrightarrow[n\rightarrow\infty]{\|.\|_X}x.
 \end{align}
Since $\varphi$ is
$\mathcal{D}$-asymptotically compact in $Y$, then there is a $y\in Y$ such that up to a subsequence,
\begin{align}\label{ff4}
 \varphi(t_n, \tau-t_n, \vartheta_{-t_n}\omega, x_n)\xrightarrow[n\rightarrow\infty]{\|.\|_Y}y.
 \end{align} It implies from (\ref{ff02}) that $y\in \mathcal{A}_Y(\tau,\omega)$. Then by (\emph{H2}), along with (\ref{ff3}) and (\ref{ff4}), we have $x=y\in \mathcal{A}_X(\tau,\omega)$ and thus
 $ \mathcal{A}_X(\tau,\omega)\subseteq  \mathcal{A}_Y(\tau,\omega)$ for every fixed $\tau\in\mathbb{R}$ and $\omega\in\Omega$. The inverse inclusion can be proved in the same
 way then we omit it here. Thus $\mathcal{A}_X=\mathcal{A}_Y$ as required.

 \emph{Step 2.}  We prove the attraction of  $\mathcal{A}_Y$
in $Y$ by a contradiction argument. Indeed,  if
there exist $\delta>0$, $x_n\in B(\tau-t_n,\vartheta_{-t_n}\omega)$ with $B\in\mathcal{D}$ and
$t_n\rightarrow +\infty$  such that
\begin{align}\label{ff5}
 \mbox{dist}_Y\Big(\varphi(t_n,\tau-t_n,\vartheta_{-t_n}\omega,x_n),
 \mathcal{A}_Y(\tau,\omega)\Big)\geq \delta.
 \end{align}
By the asymptotic compactness of $\varphi$ in $Y$, there exists $y_0\in Y$  such that up to a
subsequence,
\begin{align}\label{ff6}
 \varphi(t_n,\tau-t_n,\vartheta_{-t_n}\omega,x_n)\xrightarrow[n\rightarrow\infty]{\|.\|_Y}y_0.
 \end{align}
On the other hand, by condition (i), there exists a large time $T>0$ such that
\begin{align}\label{ff7}
 y_n=\varphi(T,\tau-t_n,\vartheta_{-t_n}\omega,x_n)&=\varphi(T,(\tau-t_n+T)-T,\vartheta_{-T}\vartheta_{-(t_n-T)}\omega,x_n)\notag\\
 &\in K(\tau-t_n+T
 ,\vartheta_{-(t_n-T)}\omega).
 \end{align}
Then by the cocycle property in Definition 2.1, along with (\ref{ff6}) and (\ref{ff7}), we infer that as $t_n\rightarrow\infty$,
\begin{align}
 \varphi(t_n,\tau-t_n,\vartheta_{-t_n}\omega,x_n)&=\varphi(t_n-T,\tau-(t_n-T),\vartheta_{-(t_n-T)}\omega,y_n)\rightarrow y_0\ \  \ \mbox{in}\ Y.\notag
 \end{align} Therefore by (\ref{ff02}), $y_0\in \mathcal{A}_Y(\tau,\omega)$. This implies that
 \begin{align}
 \mbox{dist}_Y\Big(\varphi(t_n,\tau-t_n,\vartheta_{-t_n}\omega,x_n),
 \mathcal{A}_Y(\tau,\omega)\Big)\rightarrow 0\notag
 \end{align} as $t_n\rightarrow\infty$, which is a contradiction to (\ref{ff5}).

\emph{Step 3.} It remains to prove the compactness of
$A_Y$ in $Y$.  Let
$\{y_n\}_{n=1}^\infty$  be a sequence in
$A_Y(\tau,\omega)$. By the invariance of $A_Y(\tau,\omega)$ which is proved in Step 1, we
 have
 $$
 \varphi(t,\tau-t,\vartheta_{-t}\omega, \mathcal{A}_Y(\tau-t,\vartheta_{-t}\omega))=\mathcal{A}(\tau,\omega).
 $$
 Then it follows that there is a sequence $\{z_n\}_{n=1}^\infty$ with $z_n\in \mathcal{A}_Y(\tau-t_n, \vartheta_{-t_n}\omega)$
such that for every $n\in\mathbb{Z}^+$,
 \begin{align}
 y_n=\varphi(t_n,\tau-t_n,\vartheta_{-t_n}\omega,z_n).\notag
 \end{align}
Note that $A_Y\in \mathcal{D}$. Then by the asymptotic compactness of $\varphi$ in $Y$, $\{y_n\}$  has a convergence subsequence in $Y$, \emph{i.e.},
 there is a $y_0\in Y$ such that
 $$
 \lim_{n\rightarrow\infty}y_n=y_0\ \ \mbox{in}\ Y.
 $$ But  $A_Y(\tau,\omega)$ is closed in $Y$, so $y_0\in A_Y(\tau,\omega)$.

  The uniqueness is easily followed by the attraction property of $\varphi$ and $A_Y\in \mathcal{D}$. This completes the total proofs. $\ \ \ \ \ \
\Box$\\

\textbf{Remark.}\ \   (i) It is pointed out  that the assumption \emph{(H1)} is necessary to guarantee that the closure of the set $\varphi(t,\tau-t,\vartheta_{-t}\omega,K(\tau-t,\vartheta_{-t}\omega))$ in $Y$ makes sense for all $t\in \mathbb{R}^+$, as in (\ref{ff02}).

(ii) We emphasize that the random attractor  $\mathcal{A}_Y$ in the non-initial space is completely determined by the absorbing set constructed in the initial space, without requiring the absorption in the non-initial space. This is different from the construction in \cite{Bao}.

\subsection{Upper semi-continuity of random attractors  in a non-initial space $Y$}

Assume that the assumptions (\emph{H1})-(\emph{H2}) hold.   Given
 the indexed set $I\subset \mathbb{R}$, for every $\varepsilon\in I$, we use $\mathcal{D}_\varepsilon$  to denote a
 a collection of some families of nonempty subsets of $X$. Let $\varphi_\varepsilon (\varepsilon\in I)$ be a continuous cocycle on  $X$ over $\mathbb{R}$ and $(\Omega,\mathcal{F},P,\{\vartheta_t\}_{t\in\mathbb{R}})$. We now consider the upper semi-continuous of pullback  attractors of a family of
  cocycle $\varphi_\varepsilon$ in $Y$.

Suppose first that  for every $t\in\mathbb{R}^+, \tau\in \mathbb{R}, \omega\in\Omega, \varepsilon_n,\varepsilon_0\in I$ with $\varepsilon_n\rightarrow \varepsilon_0,$ and
 $x_n,x\in X$ with $x_n\rightarrow x$, there holds
 \begin{align}\label{ff8}
\lim_{n\rightarrow \infty}\varphi_{\varepsilon_n}(t,\tau-t,\vartheta_{-t}\omega,x_n)=\varphi_{\varepsilon_0}(t,\tau-t, \vartheta_{-t}\omega, x)\ \ \ \mbox{in}\ X.
\end{align}
Suppose second that there exists a map $R_{\varepsilon_0}: \mathbb{R}\times \Omega\mapsto \mathbb{R}^+$ such that
the family
\begin{align}\label{ff9}
B_0=\{B_{0}(\tau,\omega)=\{x\in X;\|x\|_X\leq R_{\varepsilon_0}(\tau,\omega)\}:\tau\in \mathbb{R},\omega\in\Omega\}\ \mbox{belongs to}\ \mathcal{D}_{\varepsilon_0}.
\end{align}
And further for every $\varepsilon\in I$, $\varphi_\varepsilon$ has  $\mathcal{D}_\varepsilon$-pullback  attractor $\mathcal{A}_\varepsilon\in \mathcal{D}_\varepsilon$ in $X\cap Y$ and  a closed and measurable $\mathcal{D}_\varepsilon$-pullback absorbing set
$K_\varepsilon\in \mathcal{D}_\varepsilon$ in $X$ such that for every $\tau\in \mathbb{R},\omega\in\Omega$,
\begin{align}\label{ff10}
\limsup_{\varepsilon\rightarrow \varepsilon_0}\|K_\varepsilon(\tau,\omega)\|\leq R_{\varepsilon_0}(\tau,\omega),
\end{align} where $\|S\|_X=\sup_{x\in S}\|x\|_X$ for a set $S$. We finally assume that for every  $\tau\in \mathbb{R}, \omega\in\Omega$,
\begin{align}\label{ff11}
\cup_{\varepsilon\in I}\mathcal{A}_\varepsilon(\tau,\omega)\  \mbox{is precompact in}\ X, \mbox{ and}
\end{align}
\begin{align}\label{ff12}
\cup_{\varepsilon\in I}\mathcal{A}_\varepsilon(\tau,\omega)\  \mbox{is precompact in}\ Y.
\end{align}

Then we have the upper semi-continuity in $Y$.\\

\textbf{Theorem 2.7.} \emph{ If (\ref{ff8})-(\ref{ff11}) hold, then for each $\tau\in \mathbb{R}, \omega\in\Omega$,
\begin{align}\label{}
\lim_{\varepsilon\rightarrow\varepsilon_0}\mbox{dist}_X(\mathcal{A}_\varepsilon(\tau,\omega),\mathcal{A}_{\varepsilon_0}(\tau,\omega))=0.\notag
\end{align}
If further (\emph{H1})-(\emph{H2}) hold and conditions (\ref{ff8})-(\ref{ff12}) are satisfied, then for each $\tau\in \mathbb{R}, \omega\in\Omega$,
\begin{align}\label{}
\lim_{\varepsilon\rightarrow\varepsilon_0}\mbox{dist}_Y(\mathcal{A}_\varepsilon(\tau,\omega),\mathcal{A}_{\varepsilon_0}(\tau,\omega))=0.\notag
\end{align}}

\emph{Proof}\ \   If (\ref{ff8})-(\ref{ff11}) hold, the upper-continuous in $X$ is proved in \cite{Wang1}.  We only need to prove the upper semi-continuity of $\mathcal{A}_\varepsilon$ at $\varepsilon=\varepsilon_0$ in $Y$.

 Suppose that there exist $\delta>0$,  $\varepsilon_n\rightarrow\varepsilon_0$ and a sequence $\{y_n\}$ with $y_n\in \mathcal{A}_{\varepsilon_n}(\tau,\omega)$ such that for all
$n\in\mathbb{N}$,
\begin{align}\label{ff18}
\lim_{\varepsilon\rightarrow\varepsilon_0}\mbox{dist}_Y(y_n,\mathcal{A}_{\varepsilon_0}(\tau,\omega))\geq 2\delta.
\end{align}
Note that $y_n\in \mathcal{A}_{\varepsilon_n}(\tau,\omega)\subset \mathbb{A}(\tau,\omega)=\cup_{\varepsilon\in I}\mathcal{A}_\varepsilon(\tau,\omega)$. Then by (\ref{ff11}) and (\ref{ff12}) and using  (\emph{H1}), there exists a $y_0\in X\cap Y$ such that up to a
subsequence,
\begin{align}\label{ff19}
\lim_{n\rightarrow\infty}y_n=y_0\  \ \ \mbox{in}\ X\ (\mbox{resp. in}\ Y).
\end{align}
It suffices to show that $\mbox{dist}_Y(y_0, \mathcal{A}_{\varepsilon_0}(\tau,\omega))=0$.
Given  a positive sequence $\{t_m\}$ with $t_m\uparrow +\infty$ as $m\rightarrow\infty$.  For $m=1$, by the invariance of $\mathcal{A}_{\varepsilon_n}$, there exists a sequence $\{y_{1,n}\}$ with
$y_{1,n}\in \mathcal{A}_{\varepsilon_n}(\tau-t_1,\vartheta_{-t_1}\omega)$  such that
\begin{align}\label{ff20}
y_n=\varphi_{\varepsilon_n}(t_1,\tau-t_1,\vartheta_{-t_1}\omega,y_{1,n}),
\end{align} for each $n\in\mathbb{N}$.
Since $y_{1,n}\in \mathcal{A}_{\varepsilon_n}(\tau-t_1,\vartheta_{-t_1}\omega)\subset \mathbb{A}(\tau-t_1,\vartheta_{-t_1}\omega)$, then by by (\ref{ff11}) and (\ref{ff12}) and using (\emph{H2}),
there is a $z_1\in X\cap Y$ and a subsequence of $\{y_{1,n}\}$ such that
\begin{align}\label{ff21}
\lim_{n\rightarrow\infty}y_{1,n}=z_1\  \ \ \ \mbox{in}\ X\ (\mbox{resp. in}\ Y).
\end{align}
Then (\ref{ff8}) and (\ref{ff21}) together imply that
\begin{align}\label{ff22}
\lim_{n\rightarrow\infty}\varphi_{\varepsilon_n}(t_1,\tau-t_1,\vartheta_{-t_1}\omega,y_{1,n})=\varphi_{\varepsilon_0}(t_1,\tau-t_1,\vartheta_{-t_1}\omega,z_1)\  \ \ \mbox{in}\ X.
\end{align}
Thus combining (\ref{ff19}), (\ref{ff20}) and (\ref{ff22}) we get that
\begin{align}\label{ff23}
y_0=\varphi_{\varepsilon_0}(t_1,\tau-t_1,\vartheta_{-t_1}\omega,z_1).
\end{align}
Note that $K_{\varepsilon_n}$ as a $\mathcal{D}_{\varepsilon_n}$-pullback absorbing set in $X$ absorbs $\mathcal{A}_{\varepsilon_n}\in \mathcal{D}_{\varepsilon_n}$,
\emph{i.e.}, there is a $T=T(\tau,\omega,\mathcal{A}_{\varepsilon_n})$ such that for all $t\geq T$,
\begin{align}\label{ff24}
\varphi(t,\tau-t,\vartheta_{-t}\omega,\mathcal{A}_{\varepsilon_n}(\tau-t,\vartheta_{-t}\omega))\subseteq K_{\varepsilon_n}(\tau,\omega).
\end{align}
Then by the invariance of $\mathcal{A}_{\varepsilon_n}(\tau,\omega)$, it follows from (\ref{ff24}) that
\begin{align}\label{ff25}
\mathcal{A}_{\varepsilon_n}(\tau,\omega)\subseteq K_{\varepsilon_n}(\tau,\omega).
\end{align}
Since $y_{1,n}\in \mathcal{A}_{\varepsilon_n}(\tau-t_1,\vartheta_{-t_1}\omega)\subseteq K_{\varepsilon_n}(\tau-t_1,\vartheta_{-t_1}\omega)$,
then by (\ref{ff21}) and (\ref{ff10}), we get that
\begin{align}\label{ff26}
\|z_1\|_X=\limsup_{n\rightarrow\infty}\|y_{1,n}\|_X \leq \limsup_{n\rightarrow\infty}\|K_{\varepsilon_n}(\tau-t_1,\vartheta_{-t_1}\omega)\|_X\leq R_{\varepsilon_0}(\tau-t_1,\vartheta_{-t_1}\omega).
\end{align}
By an induction argument, for each $m\geq1$, there is $z_m\in X\cap Y$ such that for all $m\in\mathbb{N}$,
\begin{align}\label{ff27}
y_0=\varphi_{\varepsilon_0}(t_m,\tau-t_m,\vartheta_{-t_m}\omega,z_m).
\end{align}
and
\begin{align}\label{ff28}
\|z_m\|_X\leq R_{\varepsilon_0}(\tau-t_m,\vartheta_{-t_m}\omega).
\end{align}
Thus from (\ref{ff9}) and (\ref{ff28}), for each $m\in\mathbb{N}$,
\begin{align}\label{ff29}
z_m\in B_0(\tau-t_m,\vartheta_{-t_m}\omega).
\end{align}
Consider that the pullback  attractor $\mathcal{A}_{\varepsilon_0}$ attracts every element in $\mathcal{D}_{\varepsilon_0}$ in the topology of $Y$
 and connection with $B_0\in\mathcal{D}_{\varepsilon_0}$.
 Then $\mathcal{A}_{\varepsilon_0}$ attracts $B_0$ in the topology of $Y$. Therefore by (\ref{ff27}) and (\ref{ff29}) we have
\begin{align}
\mbox{dist}_Y(y_0, \mathcal{A}_{\varepsilon_0}(\tau,\omega))=\mbox{dist}_Y(\varphi_{\varepsilon_0}(t_m,\tau-t_m,\vartheta_{-t_m}\omega,z_m), \mathcal{A}_{\varepsilon_0}(\tau,\omega))\rightarrow0,\notag
\end{align}
as $m\rightarrow\infty$. That is to say, $\mbox{dist}_Y(y_0, \mathcal{A}_{\varepsilon_0}(\tau,\omega))=\inf_{u\in\mathcal{A}_{\varepsilon_0}(\tau,\omega)}\|y_0-u\|_Y=0$ and thus we can choose a $u_0\in\mathcal{A}_{\varepsilon_0}(\tau,\omega)$
such that
\begin{align}\label{ff30}
\|y_0-u_0\|_Y\leq \delta.
\end{align}
Therefore, by (\ref{ff19}) and (\ref{ff30}), as $n\rightarrow\infty$,
$$
\mbox{dist}_Y(y_n, \mathcal{A}_{\varepsilon_0}(\tau,\omega))\leq \|y_n-u_0\|_Y\leq  \|y_n-y_0\|_Y+\delta\rightarrow\delta,
$$ which is a contradiction to (\ref{ff18}). This concludes the proof.  $\ \ \ \ \ \ \Box$\\\\

We next consider a special case of Theorem 2.7 above, in which case  the limit cocycle $\varphi_{\varepsilon_0}$ is independent of the parameter $\omega\in\Omega$.
We call such  $\varphi_{\varepsilon_0}$ a deterministic  non-autonomous cocycle on $X$ over $\mathbb{R}$. This is, $\varphi_{\varepsilon_0}$ satisfies the following two statements:

(i) $\varphi_0(0,\tau,.)$ is the identity on $X$;

(ii) $\varphi_0(t+s,\tau,.)=\varphi_0(t,\tau+s,.)\circ \varphi_0(s,\tau,.)$.

If $\varphi_0(t, \tau, .): X\mapsto X$ is continuous, then  $\varphi_{\varepsilon_0}$ is called a  deterministic  non-autonomous  continuous cocycle on $X$ over $\mathbb{R}$.

 Let $\mathcal{D}_{\varepsilon_0}$ be a collection  of some
families of nonempty subsets of $X$ denoted by
$$
\mathcal{D}_{\varepsilon_0}=\{B=\{B(\tau)\neq\emptyset; B(\tau)\in 2^X, \tau\in\mathbb{R}\}\}.
$$
A family $\mathcal{A}_{\varepsilon_0}\in \mathcal{D}_{\varepsilon_0}$ is called a $\mathcal{D}_{\varepsilon_0}$-pullback attractor of $\varphi_{\varepsilon_0}$
in $X$ (resp. in $Y$) if

(i) for each $\tau\in\mathbb{R}$, $\mathcal{A}_{\varepsilon_0}(\tau)$ is compact in $X$(resp. of $Y$);

(ii) $\varphi_{\varepsilon_0}(t,\tau,\mathcal{A}_{\varepsilon_0}(\tau))=\mathcal{A}_{\varepsilon_0}(\tau+t)$ for all $t\in\mathbb{R}^+$ and $\tau\in\mathbb{R}$;

(iii)  $\mathcal{A}_{\varepsilon_0}$ pullback attracts every element of $\mathcal{D}_{\varepsilon_0}$ under the Hausdorff semi-metric of  $X$ (resp. of Y).

In order to obtain the convergence at $\varepsilon=\varepsilon_0$ in $Y$, we make some modifications of the conditions used in random case.
We assume that  for every $t\in\mathbb{R}^+, \tau\in \mathbb{R}, \omega\in\Omega, \varepsilon_n\in I$ with $\varepsilon_n\rightarrow \varepsilon_0,$ and
 $x_n,x\in X$ with $x_n\rightarrow x$, there
 holds
 \begin{align}\label{ff31}
\lim_{n\rightarrow \infty}\varphi_{\varepsilon_n}(t,\tau-t,\vartheta_{-t}\omega,x_n)=\varphi_{\varepsilon_0}(t,\tau-t, x)\ \ \ \mbox{in}\ X.
\end{align}
There exists a map $R^\prime_{\varepsilon_0}: \mathbb{R}\mapsto \mathbb{R}$ such that
the family
\begin{align}\label{ff32}
B^\prime_0=\{B^\prime_{0}(\tau)=\{x\in X;\|x\|_X\leq R^\prime_{\varepsilon_0}(\tau)\};\tau\in \mathbb{R}\}\ \mbox{belongs to}\ \mathcal{D}_{\varepsilon_0}.
\end{align}
For every $\varepsilon\in I$,  $\varphi_\varepsilon$ has a closed measurable $\mathcal{D}_{\varepsilon}$-pullback absorbing set
$K_\varepsilon=\{K_\varepsilon(\tau,\omega);\omega\in\Omega\}\in \mathcal{D}_\varepsilon$ in $X$ such that for every $\tau\in \mathbb{R},\omega\in\Omega$,
\begin{align}\label{ff33}
\limsup_{\varepsilon\rightarrow \varepsilon_0}\|K_\varepsilon(\tau,\omega)\|\leq R^\prime_{\varepsilon_0}(\tau).
\end{align}

Then we have the following, which can be proved by a similar argument as Theorem 2.7 and so the proof is omitted. \\

 \textbf{Theorem 2.8. } \emph{ If (\ref{ff11}) and (\ref{ff31})-(\ref{ff33}) hold, then
  for each $\tau\in \mathbb{R}, \omega\in\Omega$,
\begin{align}\label{}
\lim_{\varepsilon\rightarrow\varepsilon_0}\mbox{dist}_X(\mathcal{A}_\varepsilon(\tau,\omega),\mathcal{A}_{\varepsilon_0}(\tau))=0.\notag
\end{align}
If  further (\emph{H1})-(\emph{H2}) hold and conditions (\ref{ff12}) and  (\ref{ff31})-(\ref{ff33}) are satisfied, then  for each $\tau\in \mathbb{R}, \omega\in\Omega$,
\begin{align}\label{}
\lim_{\varepsilon\rightarrow\varepsilon_0}\mbox{dist}_Y(\mathcal{A}_\varepsilon(\tau,\omega),\mathcal{A}_{\varepsilon_0}(\tau))=0.\notag
\end{align}}

\section{ Non-autonomous reaction-diffusion equation on $\mathbb{R}^N$ with multiplicative noise}

For the non-autonomous reaction-diffusion equation (\ref{eq1})-(\ref{eq2}), the nonlinearity  $f(x,s)$
satisfies almost the same assumptions as \cite{Wang1}, \emph{i.e.},
for $ x\in\mathbb{R}^N$ and $s\in\mathbb{R}$,
\begin{align}\label {a1}
&f(x,s)s\leq -\alpha_1|s|^p+\psi_1(x),
\end{align}
\begin{align}\label {a2}
&|f(x,s)|\leq \alpha_2|s|^{p-1}+\psi_2(x),
\end{align}
\begin{align}\label {a3}
&\frac{\partial f}{\partial s}(f(x,s)\leq \alpha_3 ,
\end{align}
\begin{align}\label {a4}
&|\frac{\partial f}{\partial x}(f(x,s)|\leq \psi_3(x).
\end{align} where $\alpha_i> 0 (i=1,2,3)$ are determined constants, $p\geq2$, $\psi_1\in L^{1}(\mathbb{R}^N)\cap L^{p/2}(\mathbb{R}^N)$,
$\psi_2\in L^2(\mathbb{R}^N)$ and $\psi_3\in L^2(\mathbb{R}^N)$. And the non-autonomous term $g$  satisfies that for every $\tau\in\mathbb{R}$ and some $\delta\in [0,\lambda)$,
\begin{align}\label {a5}
\int_{-\infty}^\tau e^{\delta s} \|g(s,.)\|_{L^2(\mathbb{R}^N)}^2ds<+\infty,
\end{align} where $\lambda$ is as in (\ref{eq1}), which implies that
\begin{align}\label {a66}
\int_{-\infty}^0 e^{\delta s} \|g(s+\tau,.)\|_{L^2(\mathbb{R}^N)}^2ds<+\infty,\ \mbox{and}\  \ g\in L^2_{Loc}(\mathbb{R},L^2(\mathbb{R}^N)).
\end{align}

For the probability space $(\Omega,\mathcal{F},P)$, we write
 $\Omega=\{\omega\in
C(\mathbb{R},\mathbb{R}); \omega(0) =0\}$. Let
 $\mathcal {F}$ be the
 Borel $\sigma$-algebra induced by the compact-open topology of
$\Omega$ and ${P}$ be  the corresponding
 Wiener measure on $(\Omega,\mathcal{F})$.
  We define a shift operator $\vartheta$ on $\Omega$ by
$$
\vartheta_t\omega(s)=\omega(s+t)-\omega(t), \ \mbox{for every}\ \omega\in\Omega, t,s\in\mathbb{R}.
$$
 Then $(\Omega,\mathcal{F},P,\{\vartheta_t\}_{t\in\mathbb{R}})$ is a measurable dynamical system. By the law of the iterated logarithm (see \cite{Rand1}), we know that
\begin{align} \label{3.6}
\frac{\omega(t)}{t}\rightarrow0, \ \mbox{as}\ |t|\rightarrow+\infty.
\end{align}

For $\omega\in \Omega$, set $z(t,\omega)=z_\varepsilon(t,\omega)=e^{-\varepsilon \omega (t)}$. Then we have
$dz+\varepsilon z\circ d\omega(t)=0.$
Put $v(t,\tau,\omega, v_0)=z(t,\omega)u(t,\tau,\omega,u_0)$, where $u$ is a solution of problem
(\ref{eq1})-(\ref{eq2}) with initial $u_0$. Then $v$ solves the follow non-autonomous equation
\begin{align} \label{pr1}
\frac{dv}{dt}+\lambda v-\Delta v=z(t,\omega)f(x,z^{-1}(t,\omega)v)+z(t,\omega)g(t,x),
\end{align}
with initial condition
\begin{align} \label{pr2}
 v(\tau,x)=v_0(x)=z(\tau,\omega)u_0(x).
\end{align}

As pointed out in \cite{Wang1}, for every $v_0\in L^2(\mathbb{R}^N)$ we may show that the problem (\ref{pr1})-(\ref{pr2})  possesses a continuous  solution $v(.)$ on $L^2(\mathbb{R}^N)$ such that
$v(.)\in C([\tau,+\infty),L^2(\mathbb{R}^N))\cap L^2_{loc}((\tau,+\infty), H^1(\mathbb{R}^N))\cap L^p_{loc}((\tau,+\infty), L^p(\mathbb{R}^N)).$
In addition, the solution $v$ is $(\mathcal{F},\mathcal{B}(L^2(\mathbb{R}^N)))$-measurable in $\Omega$. Then formally $u(.)=z^{-1}(.,\omega)v(.)$ is a $(\mathcal{F},\mathcal{B}(L^2(\mathbb{R}^N)))$-measurable and continuous solution of problem (\ref{eq1})-(\ref{eq2}) on $L^2(\mathbb{R}^N)$  with $u_0=z^{-1}(\tau,\omega)v_0$.

Define
the mapping $\varphi: \mathbb{R}^+\times\mathbb{R}\times \Omega \times
L^2(\mathbb{R}^N)\rightarrow L^2(\mathbb{R}^N)$  such that
\begin{align}\label{eq0}
 \varphi(t,\tau,\omega, u_0)&=u(t+\tau,\tau,\vartheta_{-\tau}\omega, u_0)=z^{-1}(t+\tau,\vartheta_{-\tau}\omega)v(t+\tau,\tau,\vartheta_{-\tau}\omega, z(\tau,\vartheta_{-\tau}\omega)u_0),
\end{align}for $u_0\in L^2(\mathbb{R}^N)$ and  $t\in\mathbb{R}^+, \tau\in\mathbb{R}, \omega\in \Omega$. Then by the measurability and
 continuity of $v$ in  $v_0\in L^2(\mathbb{R}^N)$ and $t\in \mathbb{R}^+$, we see
 that the mappings $\varphi$  is $(\mathcal{B}(\mathbb{R}^+)\times \mathcal{F}\times \mathcal{B}(L^2(\mathbb{R}^N)))\mapsto \mathcal{B}(L^2(\mathbb{R}^N))$-measurable. That is to say, the mappings $\varphi$ defined by (\ref{eq0}) is
 a continuous cocycle on $L^2(\mathbb{R}^N)$  over $\mathbb{R}$ and $(\Omega,\mathcal{F},P,\{\vartheta_t\}_{t\in\mathbb{R}})$.
 Furthermore, from (\ref{eq0}) we infer that
\begin{align} \label{eq00}
 \varphi(t,\tau-t,\vartheta_{-t}\omega, u_0)&=u(\tau,\tau-t,\vartheta_{-\tau}\omega, u_0)=z(-\tau,\omega)v(\tau,\tau-t,\vartheta_{-\tau}\omega, z(\tau-t,\vartheta_{-\tau}\omega)u_0).
\end{align}

 We define the collection  $\mathcal{D}$ as
\begin{align} \label{D}
 &\mathcal{D}=\{B=\{B(\tau,\omega)\subseteq L^2(\mathbb{R}^N); \tau\in\mathbb{R}, \omega\in \Omega\};\notag\\&\ \  \ \ \ \ \ \ \ \lim\limits_{t\rightarrow+\infty}e^{-\delta t}\|B(\tau-t,\vartheta_{-t}\omega)\|^2=0\ \mbox{for}\ \tau\in\mathbb{R}, \omega\in \Omega,\delta<\lambda\}
\end{align} where $\|B\|=\sup_{v\in B}\|v\|_{L^2(\mathbb{R}^N)}$ and $\lambda$ is in (\ref{pr1}).
Note that this collection $\mathcal{D}$ is much larger that
the collection defined by \cite{Wang1}. That  is to say, the collection $\mathcal{D}$ defined above includes all tempered families of bounded nonempty
subsets of $L^2(\mathbb{R}^N)$.

We can show that all the results in \cite{Wang1} hold for this collection $\mathcal{D}$ defined by (\ref{D}).  Thus, the existence and upper semi-continuous of
$\mathcal{D}$-pullback  attractors for the cocycle $\varphi_\varepsilon$ in \emph{the initial space} $L^2(\mathbb{R}^N)$ have been proved by \cite{Wang1}.\\

\textbf{Theorem 3.1(\cite{Wang1})}.\emph{ Assume that (\ref{a1})-(\ref{a5}) hold. Then the cocycle $\varphi_\varepsilon$ has a unique $\mathcal{D}$-pullback  attractor $\mathcal{A}_\varepsilon=\{\mathcal{A}_\varepsilon(\tau,\omega),\tau\in\mathbb{R}, \omega\in \Omega\}$ in $L^2(\mathbb{R}^N)$, given by
\begin{align}\label{L2}
\mathcal{A}_{\varepsilon}(\tau,\omega)=\bigcap_{s>0}\overline{\bigcup_{t\geq s} \varphi(t,\tau-t,\vartheta_{-t}\omega, K_\varepsilon(\tau-t,\vartheta_{-t}\omega))}^{L^2(\mathbb{R}^N)}, \ \
\tau\in \mathbb{R},\omega\in \Omega,
\end{align} where $K_\varepsilon$ is a closed and measurable $\mathcal{D}$-pullback absorbing set of $\varphi_\varepsilon$ in $L^2(\mathbb{R}^N)$. Furthermore,
$\mathcal{A}_\varepsilon$ is upper semi-continuous in $L^2(\mathbb{R}^N)$ at $\varepsilon=0$.
}\\

 Note that in  most cases, we  write $v$ (resp. $\varphi$ and $z$) as the abbreviation of $v_\varepsilon$(resp. $\varphi_\varepsilon$ and $z_\varepsilon$).

 In the following, we consider the applications of Theorem 2.6-2.8 to the non-autonomous stochastic reaction-diffusion (\ref{eq1})-(\ref{eq2}). We will strengthen the result of
 Theorem 3.1 holds in  the smooth functions space $H^1(\mathbb{R}^N)$. In particular, we prove  the upper semi-continuity of the obtained attractors
 $\mathcal{A}_\varepsilon$ in $H^1(\mathbb{R}^N)$.

\section{Existence of pullback attractor in $H^{1}(\mathbb{R}^N)$}

In this section, we apply Theorem 2.6 to  prove the existence  of $\mathcal{D}$-pullback  attractors in $H^1(\mathbb{R}^N)$ for the cocycle defined in (\ref{eq0}).
 To this end,
we need to prove the uniform smallness of solutions outside a large ball under  $H^1(\mathbb{R}^N)$ norm (see Lemma 4.4), and in the bounded ball of $\mathbb{R}^N$, we will prove the asymptotic compactness of solutions by a combined space splitting  and function truncation technique (see Lemma 4.5 and Lemma 4.6).

Consider that
$e^{-|\omega(s)|}\leq z(s,\omega)=e^{-\varepsilon \omega(s)}\leq e^{|\omega(s)|}$ for $\varepsilon\in(0,1]$, and $\omega(s)$ is continuous function in $s$. Then there exist two positive random constants
 $E=E(\omega)$ and  $F=F(\omega)$ depending only on $\omega$ such that for all $s\in[-1,0]$ and $\varepsilon\in(0,1]$.
\begin{align} \label{EE}
&0<E\leq z(s,\omega)\leq F,\ \,\omega\in\Omega.
\end{align}

Hereafter, we denote by $\|.\|, \|.\|_p$ and $\|.\|_{H^1}$ the norms in $L^2(\mathbb{R}^N),L^p(\mathbb{R}^N)$ and $H^1(\mathbb{R}^N)$, respectively.
The number $c$  is a generic positive constant  independent of $\tau,\omega,B$ and $\varepsilon$ in any place. We always assume $p>2$ in the following discussions.

\subsection{$H^1$-tail estimate of solutions}

This can be achieved by a series of previously proved lemmas. First we stress that Lemma 5.1 in \cite{Wang1} holds on the compact interval $[\tau-1,\tau]$, which is necessary for
us to estimate of the tail of solutions in $H^1(\mathbb{R}^N)$.\\

\textbf{Lemma 4.1.}
 \emph{Assume that (\ref{a1})-(\ref{a5}) hold. Given  $\tau\in\mathbb{R}, \omega\in\Omega$ and  $B=\{B(\tau,\omega);\tau\in\mathbb{R},\omega\in\Omega\}\in\mathcal{D}$, then there exists a constant  $T=T(\tau,\omega,B)\geq2$
such that for all $t\geq T$,  the solution $v$ of problem (\ref{pr1})-(\ref{pr2}) satisfies that for every $\xi\in[\tau-1,\tau]$,
\begin{align} \label{ee1}
\|v(\xi,\tau-t,\vartheta_{-\tau}\omega, &z(\tau-t,\vartheta_{-\tau}\omega)u_0)\|_{H^1(\mathbb{R}^N)}^2\leq L_1(\tau,\omega,\varepsilon),
\end{align}
\begin{align} \label{ee2}
\int_{\tau-t}^\tau e^{\lambda(s-\tau)}\Big(\|v(s,\tau-t,\vartheta_{-\tau}\omega, v_0)\|_{H_1}^2+z^{2-p}(s,\omega)\|v(s,\tau-t,\vartheta_{-\tau}\omega, v_0)\|_p^p\Big)ds\leq  L_1(\tau,\omega,\varepsilon),
\end{align}
where $u_{0}\in B(\tau-t,\vartheta_{-t}\omega)$ and $L_1(\tau,\omega,\varepsilon)=cz^{-2}(-\tau,\omega)\int_{-\infty}^0 e^{\lambda s}z^2(s,\omega)( \|g(s+\tau,.)\|^2+1)ds$.}\\

\emph{Proof}\ \  By (\ref{pr1}), it is easy to calculate that
\begin{align} \label{se01}
\frac{d}{dt}\|v\|^2+\frac{3}{2}\lambda\|v\|^2+\|\nabla v\|^2+\alpha_1z^{2-p}(t,\omega)\|v\|_p^p\leq c z^2(t,\omega)(\|g(t,.)\|^2+\|\psi_1\|_1),
\end{align}
to which we apply Gronwall'lemma over the interval $[\tau-t,\xi]$, where $\xi\in[\tau-1,\tau]$ and $t\geq2$, we find that, along with  $\omega$ replaced by $\vartheta_{-\tau}\omega$,
\begin{align} \label{se02}
&\|v(\xi,\tau-t,\vartheta_{-\tau}\omega, v_0)\|^2+\frac{\lambda}{2}\int_{\tau-t}^\xi e^{\lambda(s-\xi)}\|v(s,\tau-t,\vartheta_{-\tau}\omega, v_0)\|^2ds\notag\\&+\int_{\tau-t}^\xi e^{\lambda(s-\xi)}\Big(\|\nabla v(s,\tau-t,\vartheta_{-\tau}\omega, v_0)\|^2
+\alpha_1z^{2-p}(s,\vartheta_{-\tau}\omega)\|v(s,\tau-t,\vartheta_{-\tau}\omega, v_0)\|_p^p\Big)ds\notag\\
&\leq e^{-\lambda(\xi-\tau+t)}z^2(\tau-t,\vartheta_{-\tau}\omega)\|u_0\|^2
+c\int_{\tau-t}^\xi e^{-\lambda(\xi-\sigma)}z^2(\sigma,\vartheta_{-\tau}\omega)(\|g(\sigma,.)\|^2+\|\psi_1\|_1)d\sigma.
\end{align}
If $\xi\in[\tau-1,\tau]$, then
\begin{align} \label{xx}
e^{-\lambda\tau}\leq e^{-\lambda\xi}\leq e^{-\lambda(\tau-1)},
\end{align}
 and therefore (\ref{se02}) along with (\ref{xx}) implies that
\begin{align} \label{se022}
&\|v(\xi,\tau-t,\vartheta_{-\tau}\omega, v_0)\|^2\notag\\&+\int_{\tau-t}^\xi e^{\lambda(s-\tau)}\Big(h\| v(s,\tau-t,\vartheta_{-\tau}\omega, v_0)\|_{H^1}^2
+\alpha_1z^{2-p}(s,\vartheta_{-\tau}\omega)\|v(s,\tau-t,\vartheta_{-\tau}\omega, v_0)\|_p^p\Big)ds\notag\\
&\leq e^{\lambda}z^2(\tau-t,\vartheta_{-\tau}\omega)\|u_0\|^2
+ce^{\lambda}\int_{\tau-t}^\tau e^{-\lambda(\tau-\sigma)}z^2(\sigma,\vartheta_{-\tau}\omega)(\|g(\sigma,.)\|^2+\|\psi_1\|_1)d\sigma\notag\\
&\leq e^{\lambda} z^{-2}(-\tau,\omega)\Big(e^{-\lambda t}z^2(-t,\omega)\|u_0\|^2
+c\int_{-\infty}^0 e^{\lambda \sigma}z^2(\sigma,\omega)( \|g(\sigma+\tau,.)\|^2+1)d\sigma\Big),
\end{align} where $h=\min\{\frac{\lambda}{2},1\}$.
From (\ref{3.6}), we can calculate that
\begin{align} \label{jjj}
\lim\limits_{t\rightarrow+\infty}z^2(-t,\omega)e^{-(\lambda-\delta)t}=0
\end{align} for $\lambda>\delta>0$.
Then  for each $\tau\in\mathbb{R},\omega\in\Omega$ and $u_0\in B(\tau-t,\vartheta_{-t}\omega)$, by (\ref{D})  we deduce that
\begin{align} \label{se03}
\lim\limits_{t\rightarrow+\infty}z^2(-t,\omega)e^{-\lambda t}\|u_0\|^2=\lim\limits_{t\rightarrow+\infty}e^{-\delta t}\|u_0\|^2=0.
\end{align}
Then by  (\ref{se022})-(\ref{se03}), it implies that there exists $T_1=T_1(\tau,\omega,B)\geq2$ such that  for all $t\geq T_1$ and $\xi\in[\tau-1,\tau]$,
\begin{align} \label{se04}
&\|v(\xi,\tau-t,\vartheta_{-\tau}\omega, v_0)\|^2\leq c z^{-2}(-\tau,\omega)\int_{-\infty}^0 e^{\lambda s}z^2(s,\omega)( \|g(s+\tau,.)\|^2+1)ds,
\end{align}
and
\begin{align} \label{se044}
&\int_{\tau-t}^\xi  e^{\lambda(s-\tau)}\Big(\|v(s,\tau-t,\vartheta_{-\tau}\omega, v_0)\|_{H_1}^2+z^{2-p}(s,\omega)\|v(s,\tau-t,\vartheta_{-\tau}\omega, v_0)\|_p^p\Big)ds\notag\\& \ \ \ \ \ \ \ \ \ \ \ \ \ \ \ \ \ \ \leq c z^{-2}(-\tau,\omega)\int_{-\infty}^0 e^{\lambda s}z^2(s,\omega)( \|g(s+\tau,.)\|^2+1)ds.
\end{align}
On the other hand, by (\ref{pr1}), we  deduce that
\begin{align} \label{H11}
\frac{d}{dt}\|\nabla v\|^2+\lambda\|\nabla v\|^2\leq  c\|\nabla v\|^2+z^2(t,\omega)(\|g(t,.)\|^2+\|\psi_3\|^2).
\end{align}
Note that $\xi-\tau+t\geq t-1\geq 1$ for $\xi\in[\tau-1,\tau]$ and $t\geq2$. Then applying Lemma 5.1 in \cite{Zhao0} over the interval $[\tau-t, \xi] $ for $\xi\in[\tau-1,\tau]$, we get that
\begin{align} \label{ttt}
\|\nabla &v(\xi,\tau-t,\vartheta_{-\tau}\omega, v_0)\|^2\notag\\
&\leq \frac{e^\lambda}{\xi-\tau+t}\int_{\tau-t}^\tau e^{\lambda(\sigma-\tau)}\|\nabla v(\sigma,\tau-t,\vartheta_{-\tau}\omega,v_0)\|^2d\sigma\notag\\
&\ \ \ \ \ +c\int_{\tau-t}^\tau e^{\lambda(\sigma-\tau)}\|\nabla v(\sigma,\tau-t,\vartheta_{-\tau}\omega, v_0)\|^2d\sigma\notag\\
&\ \ \ \ \ \ \ \ +c\int_{\tau-t}^\tau e^{\lambda(\sigma-\tau)}z^{2}(\sigma,\vartheta_{-\tau}\omega)(\|g(\sigma,.)\|^2+1)d\sigma\notag\\
&\ \ \ \ \ \ \ \  \ \ \ \leq c\int_{\tau-t}^\tau e^{\lambda(\sigma-\tau)}\|\nabla v(\sigma,\tau-t,\vartheta_{-\tau}\omega, v_0)\|^2d\sigma\notag\\
&\ \ \ \ \ \ \ \  \ \ \ \ \ \  +c z^{-2}(-\tau,\omega)\int_{-\infty}^0 e^{\lambda s}z^2(s,\omega)( \|g(s+\tau,.)\|^2+1)ds,
\end{align} where we  have used (\ref{xx}). Then by (\ref{ttt}) and (\ref{se044}) it follows that  there exists $T_2=T_2(\tau,\omega,B)\geq2$ such that for all $t\geq T_2$,
\begin{align} \label{qw}
&\|\nabla v(\xi,\tau-t,\vartheta_{-\tau}\omega, v_0)\|^2 \leq c z^{-2}(-\tau,\omega)\int_{-\infty}^0 e^{\lambda s}z^2(s,\omega)( \|g(s+\tau,.)\|^2+1)ds,
\end{align}
which is finite for all $\xi\in[\tau-1,\tau]$. Taking $T=\max\{T_1,T_2\}$, then for all $t\geq T$, (\ref{se04}) and (\ref{qw}) together imply the desired.
 $\ \ \ \ \ \ \Box$\\

\textbf{Lemma 4.2.} \emph{Assume that (\ref{a1})-(\ref{a5}) hold. Given  $\tau\in\mathbb{R}, \omega\in\Omega$ and  $B=\{B(\tau,\omega);\tau\in\mathbb{R},\omega\in\Omega\}\in\mathcal{D}$, then for every $\epsilon>0$, there exist two constants $T=T(\tau,\omega,\epsilon,B)\geq2$ and $R=R(\tau,\omega,\epsilon)>1$
such that the weak solution $v$ of problem (\ref{pr1})-(\ref{pr2}) satisfies that for all $t\geq T$ and $k\geq R$,
\begin{align}\label{}
\int\limits_{|x|\geq k}| &v(\tau,\tau-t,\vartheta_{-\tau}\omega, z(\tau-t,\vartheta_{-\tau}\omega)u_0)|^2dx\notag\\
&+\int_{\tau-t}^\tau e^{\lambda(s-\tau)}\int\limits_{|x|\geq k}|\nabla v(s,\tau-t,\vartheta_{-\tau}\omega, z(\tau-t,\vartheta_{-\tau}\omega)u_0)|^2dxds\leq\epsilon,\notag
\end{align}
where $u_0\in B(\tau-t,\vartheta_{-t}\omega)$,  $R$ and $T$ are independent of $\varepsilon$.}\\

\emph{Proof} \ The proof is a simple modification of  the proof of Lemma 5.5 in \cite{Wang1}.  We first need to define a smooth function  $\xi(.)$ on $\mathbb{R}^+$ such that
\begin{equation} \xi(s)=
\begin{cases}0,& \mbox{if}\ 0\leq s\leq 1,\\
0\leq\xi(s)\leq1, &\mbox{if}\  1\leq s\leq 2,\\
1,&\mbox{if}\ s\geq 2,
\end{cases}
\end{equation}which obviously implies that there is a  positive constant $C_1$ such that the $|\xi^\prime(s)|+|\xi^{\prime\prime}(s)|\leq C_1$ for all $s\geq0$.
For convenience, we write $\xi=\xi(\frac{|x|^2}{k^2})$.

From (\ref{pr1}), we know that
\begin{align}\label{4.3}
\frac{1}{2}\frac{d}{dt}\int\limits_{\mathbb{R}^N}\xi|v|^2dx&+\lambda\int\limits_{\mathbb{R}^N}\xi|v|^2dx-\int\limits_{\mathbb{R}^N}\xi v\Delta v dx
\notag\\&=z(t,\omega)\int\limits_{\mathbb{R}^N}f(x,u)v\xi dx+z(t,\omega)\int\limits_{\mathbb{R}^N}gv\xi dx.
\end{align}
By calculation, we have the following:
\begin{align}\label{4.301}
&\ \  \ \ \ \ \ \ \ \ \ -\int\limits_{\mathbb{R}^N}\xi v\Delta v dx=\int\limits_{\mathbb{R}^N}v(\nabla\xi.\nabla v)dx+\int\limits_{\mathbb{R}^N}\xi|\nabla v|^2dx\notag\\
&\geq -\Big|\int\limits_{k\leq|x|\leq\sqrt{2}k}\xi^\prime\frac{|x|}{k^2}|v||\nabla v|dx\Big|+\int\limits_{\mathbb{R}^N}\xi|\nabla v|^2dx\geq-\frac{c}{k}\|v\|^2_{H^1}+\int\limits_{\mathbb{R}^N}\xi|\nabla v|^2dx,
\end{align}
\begin{align}\label{4.302}
z(t,\omega)\int\limits_{\mathbb{R}^N}f(x,u)v\xi dx\leq -\alpha_1z^{2-p}(t,\omega)\int\limits_{\mathbb{R}^N}\xi|v|^p dx+z^{2}(t,\omega)\int\limits_{\mathbb{R}^N}\xi\psi_1 dx,
\end{align}
\begin{align}\label{4.303}
\Big|z(t,\omega)\int\limits_{\mathbb{R}^N}gv\xi dx\Big|\leq \frac{\lambda}{2}\int\limits_{\mathbb{R}^N}\xi|v|^2dx+\frac{1}{2\lambda} z^{2}(t,\omega)\int\limits_{\mathbb{R}^N}\xi g^2dx.
\end{align}
Then combining (\ref{4.3})- (\ref{4.303}) we get that
\begin{align}\label{4.4}
\frac{d}{dt}\int\limits_{\mathbb{R}^N}\xi|v|^2dx&+\lambda\int\limits_{\mathbb{R}^N}\xi|v|^2dx+\int\limits_{\mathbb{R}^N}\xi |\nabla v|^2 dx
\notag\\&\ \ \ \  \ \ \ \ \ \leq cz^2(t,\omega)\int\limits_{\mathbb{R}^N}\xi(|\psi_1|+|g(t,x)|^2)dx+\frac{c}{k}\|v\|^2_{H^1}.
\end{align}
Applying the Gronwall'lemma to (\ref{4.4}) over $[\tau-t,\tau]$,  we find that,  along with  $\omega$ replaced by $\vartheta_{-\tau}\omega$,
\begin{align}\label{4.5}
\int\limits_{\mathbb{R}^N}\xi|v(\tau,\tau-t, &\vartheta_{-\tau}\omega, v_0)|^2dx
+\int_{\tau-t}^\tau e^{\lambda(s-\tau)} \int\limits_{\mathbb{R}^N}\xi|\nabla v(\tau,\tau-t,\vartheta_{-\tau}\omega,v_0)|^2dxds\notag\\
&\leq cz^{-2}(\tau,\omega)\int_{-\infty}^0e^{\lambda s} z^2(s,\omega)\int\limits_{|x|\geq k}(|\psi_1|+|g(s+\tau,x)|^2)dxds\notag\\
&+\frac{c}{k}\int_{\tau-t}^\tau e^{\lambda(\tau-s)}\|v(s,\tau-t,\vartheta_{-\tau}\omega, v_0)\|^2_{H^1}ds+z^2(-\tau,\omega) e^{-\lambda t}z^2(-t,\omega)\|u_0\|^2.
\end{align}
According to Lemma 4.1,
there exist $T_1=T_1(\tau,\omega,B)\geq2$ and $R_1=R_1(\tau,\omega,\epsilon)>2$ such that for all
 $t\geq T_1$  and $k\geq R_1$,
\begin{align}\label{4.8}
\frac{c}{k}\int_{\tau-t}^\tau e^{\lambda(\tau-s)}\|v(s,\tau-t,\vartheta_{-\tau}\omega, v_0)\|^2_{H^1}ds\leq \frac{cL(\tau,\omega,\varepsilon)}{k}\leq \frac{\epsilon}{3}.
\end{align}
On the other hand, for each $\tau\in\mathbb{R},\omega\in\Omega$ and $u_0\in B(\tau-t,\vartheta_{-t}\omega)$,  by (\ref{se03}),
there exists $T_2=T_2(\tau,\omega,B,\epsilon)>0$ such that for all $t\geq T_2$,
\begin{align} \label{4.9}
z^2(-\tau,\omega) e^{-\lambda t}z^2(-t,\omega)\|u_0\|^2\leq \frac{\epsilon}{3}.
\end{align}
By (\ref{jjj}), there exists a random variable $a(\omega)$ depending only on $\omega$ such that
$$
 0<e^{(\lambda-\delta) s}z^2(s,\omega)\leq a(\omega), \ \ \mbox{for}\ s\in(-\infty,0].
$$
Then  by (\ref{a66}), we can deduce that for every $\tau\in\mathbb{R}$,
\begin{align}\label{4.10}
\int_{-\infty}^0e^{\lambda s} z^2(s,\omega)\int\limits_{\mathbb{R}^N}|g(s+\tau,x)|^2dxds
&=\int_{-\infty}^0 e^{(\lambda-\delta) s}z^2(s,\omega)e^{\delta s}\|g(s+\tau,.)\|^2ds\notag\\
&\leq a(\omega)\int_{-\infty}^0 e^{\delta s}\|g(s+\tau,.)\|^2ds<+\infty,
\end{align} where $\delta\in[0,\lambda)$.
Then by (\ref{4.10}) and  $\psi_1\in L^1$,  there exists $R_2=R_2(\tau,\omega,\epsilon)$ such that for all $k\geq R_2$,
\begin{align}\label{4.11}
cz^{-2}(\tau,\omega)\int_{-\infty}^0e^{\lambda s} z^2(s,\omega)\int\limits_{|x|\geq k}(|\psi_1|+|g(s+\tau,x)|^2)dxds\leq\frac{\epsilon}{3}.
\end{align}
Given $T=\max\{T_1,T_2\}$ and ${R}=\max\{R_1,R_2\}$, then combining (\ref{4.8})-(\ref{4.9})  and (\ref{4.11}) into (\ref{4.5}), we have for all $t\geq T$ and $k\geq {R}$,
\begin{align}\label{}
\int\limits_{|x|\geq k}&|v(\tau,\tau-t, \vartheta_{-\tau}\omega, v_0)|^2dx+\int_{\tau-t}^\tau e^{\lambda(s-\tau)} \int\limits_{|x|\geq k}|\nabla v(\tau,\tau-t,\vartheta_{-\tau}\omega,v_0)|^2dxds
\leq \epsilon.\notag
\end{align}
Then the desired result follows. $\ \ \ \ \ \ \Box$\\

\textbf{Lemma 4.3.} \emph{Assume that (\ref{a1})-(\ref{a5}) hold. Given  $\tau\in\mathbb{R}, \omega\in\Omega$ and  $B=\{B(\tau,\omega);\tau\in\mathbb{R},\omega\in\Omega\}\in\mathcal{D}$, then  there exists $T=T(\tau,\omega,B)\geq 2$
such that the weak solution $v$ of problem (\ref{pr1})-(\ref{pr2}) satisfies that for all $t\geq T$,
\begin{align}
&\int_{\tau-1}^\tau e^{\lambda(s-\tau)}z^{4-2p}(s,\vartheta_{-\tau}\omega)\| v(s,\tau-t,\vartheta_{-\tau}\omega,z(\tau-t,\vartheta_{-\tau}\omega)u_0)\|_{2p-2}^{2p-2}ds\leq L_2(\tau,\omega,\varepsilon),\notag\\
&\int_{\tau-1}^\tau e^{\lambda(s-\tau)} \| v_s(s,\tau-t,\vartheta_{-\tau}\omega, z(\tau-t,\vartheta_{-\tau}\omega)u_0)\|^2\Big)ds\leq L_3(\tau,\omega,\varepsilon),\notag
\end{align} where $v_s=\frac{\partial v}{\partial s}$, $u_0\in B(\tau-t,\vartheta_{-t}\omega)$ and
\begin{align}
&L_2(\tau,\omega,\varepsilon)=cz^{-2}(-\tau,\omega)F^{2-p}\Big(b(\omega)\int_{-\infty}^{0}e^{\lambda s}z^2(s,\omega)(\|g(s+\tau,.)\|^2+1)ds\notag\\
&\ \ \ \ \ \ \ \ \ \ \ \ \ \ \ +\int_{-\infty}^{0} e^{\lambda s} z^p(s, \omega)(\|g(s+\tau,.)\|^2+1)ds\Big),\notag\\
&L_3(\tau,\omega,\varepsilon)=c(F^{2-p}b(\omega)+1) L_1(\tau,\omega,\varepsilon)+cz^{-2}(-\tau,\omega)F^{2-p}\int_{-\infty}^{0} e^{\lambda s} z^p(s, \omega)(\|g(s+\tau,.)\|^2+1)ds\Big),\notag
\end{align} where $F$ is as in (\ref{EE}), $L_1(\tau,\omega,\varepsilon)$ is as in (\ref{ee1}) and $b(\omega)$ is as in (\ref{4.1667}}).
\\

\emph{Proof}\ \
 We multiply (\ref{pr1}) by $|v|^{p-2}v$  and then integrate over $\mathbb{R}^N$ to yield  that
\begin{align} \label{4.12}
\frac{1}{p}\frac{d}{dt}\|v\|_{p}^p +\lambda\|v\|_p^p\leq z(t,\omega)\int\limits_{\mathbb{R}^N}f(x,z^{-1}v)
|v|^{p-2}vdx+z(t,\omega)\int\limits_{\mathbb{R}^N}|v|^{p-2}v{g}dx.
\end{align}
By using (\ref{a1}),  we see that
\begin{align} \label{4.13}
z(t,\omega)\int\limits_{\mathbb{R}^N}f(x,z^{-1}v) |v|^{p-2}vdx
&\leq -\alpha_1z^{2-p}(t,\omega)\int\limits_{\mathbb{R}^N}|v|^{2p-2}dx+z^2(t,\omega)\int\limits_{\mathbb{R}^N}\psi_1(x)|v|^{p-2}dx\notag\\
& \leq-\alpha_1z^{2-p}(t,\omega)\int\limits_{\mathbb{R}^N}|v|^{2p-2}dx+\frac{\lambda}{2}\|v\|_p^p+z^{p}(t,\omega)\|\psi_1\|^{p/2}_{p/2}.
\end{align}
At the same time, the last term on the right hand side  of (\ref{4.12}) is bounded by
\begin{align} \label{4.14}
z(t,\omega)\int\limits_{\mathbb{R}^N}|v|^{p-2}v{g}dx\leq
\frac{1}{2}\alpha_1z^{2-p}(t,\omega)\int\limits_{\mathbb{R}^N}|v|^{2p-2}dx+cz^p(t,\omega)\|g(t,.)\|^2.
\end{align}
 Combination (\ref{4.12})-(\ref{4.14}),  we obtain that
\begin{align} \label{4.15}
&\frac{d}{dt}\|v\|_{p}^p+2\lambda \|v\|_{p}^p+\alpha_1z^{2-p}(t,\omega)\|v\|_{2p-2}^{2p-2}\leq c
z^p(t,\omega)(\|g(t,.)\|^2+1).
\end{align}
 Applying Lemma 5.1 in \cite{Zhao0} over $[\tau-t,\xi]$ for $\xi\in[\tau-1,\tau]$ and $t\geq2$, along with  $\omega$ replaced by $\vartheta_{-\tau}\omega$, we deduce that
\begin{align} \label{4.16}
\|v(\xi, \tau-t,\vartheta_{-\tau}\omega, v_0)\|_{p}^p
&\leq c\int_{\tau-t}^{\tau} e^{2\lambda (s-\tau)} \|v(s,\tau-t,\vartheta_{-\tau}\omega,v_0)\|_{p}^pds\notag\\
&+cz^{-p}(-\tau,\omega)\int_{-\infty}^{0} e^{\lambda s} z^p(s, \omega)(\|g(s+\tau,.)\|^2+1)ds,
\end{align} where we  have used (\ref{xx}) and $\xi-\tau+t\geq t-1\geq 1$ for $\xi\in[\tau-1,\tau]$ and $t\geq2$.
On the other hand, we see that
\begin{align} \label{4.166}
&\ \ \ \ \ \ \ \  \int_{\tau-t}^\tau e^{2\lambda(s-\tau)}\|v(s, \tau-t, \vartheta_{-\tau}\omega, v_0)\|_{p}^pds\notag\\
&= z^{2-p}(-\tau,\omega)\int_{\tau-t}^\tau e^{2\lambda(s-\tau)}z^{p-2}(s-\tau,\omega)z^{2-p}(s,\vartheta_{-\tau}\omega)
\|v(s, \tau-t, \vartheta_{-\tau}\omega, v_0)\|_{p}^pds.
\end{align}
Consider that when $s\rightarrow-\infty$, $e^{\lambda s}z^{p-2}(s,\omega)$ and $e^{\lambda s}z^{2-p}(s,\omega)\rightarrow0$. Then there exists a variable $b(\omega)$ depending only on $\omega$ such that
\begin{align} \label{4.1667}
0<e^{\lambda s}z^{p-2}(s,\omega)+ e^{\lambda s}z^{2-p}(s,\omega)\leq b(\omega),\ \ \  s\in(-\infty,0].
\end{align}
from which  and (\ref{4.166}), association with  Lemma 4.1,  it follows that there exists $T=T(\tau,\omega,B)\geq2$ such that for all $t\geq T$,
\begin{align} \label{yy}
\int_{\tau-t}^\tau e^{2\lambda(s-\tau)}&\|v(s, \tau-t, \vartheta_{-\tau}\omega, v_0)\|_{p}^pds\notag\\
&\leq b(\omega)z^{2-p}(-\tau,\omega)\int_{\tau-t}^\tau e^{\lambda(s-\tau)}z^{2-p}(s,\vartheta_{-\tau}\omega)
\|v(s, \tau-t, \vartheta_{-\tau}\omega, v_0)\|_{p}^pds
\notag\\
&\leq cb(\omega)z^{-p}(-\tau,\omega)\int_{-\infty}^{0}e^{\lambda s}z^2(s,\omega)(\|g(s+\tau,.)\|^2+1)ds.
\end{align}
Then by (\ref{4.16}) and (\ref{yy}) we get that for all $t\geq T$ and $\xi\in [\tau-1,\tau]$,
\begin{align} \label{4.17}
\|v(\xi, \tau-t, \vartheta_{-\tau}\omega, v_0)\|_{p}^p
&\leq c z^{-p}(-\tau,\omega)\Big(b(\omega)\int_{-\infty}^{0}e^{\lambda s}z^2(s,\omega)(\|g(s+\tau,.)\|^2+1)ds\notag\\
& +\int_{-\infty}^{0} e^{\lambda s} z^p(s, \omega)(\|g(s+\tau,.)\|^2+1)ds\Big).
\end{align}
In (\ref{4.15}), omitting the number 2 of the second term on the left hand side,  we multiply (\ref{4.15})  by $e^{\lambda(t-\tau)}$  and then integrate (w.r.t $t$) from $[\tau-1,\tau]$ to yield  that, along with
$\omega$  replaced by $\vartheta_{-\tau}\omega$,
\begin{align} \label{4.20}
\int_{\tau-1}^\tau e^{\lambda(s-\tau)}z^{2-p}(s,\vartheta_{-\tau}\omega)\|v(s,\tau-t,&\vartheta_{-\tau}\omega, v_0)\|_{2p-2}^{2p-2}ds
\leq e^{-\lambda}\|v(\tau-1, \tau-t,\vartheta_{-\tau}\omega, v_0)\|_{p}^p\notag\\
&+c\int_{\tau-1}^\tau e^{\lambda(s-\tau)}z^p(s,\vartheta_{-\tau}\omega)(\|{g}(s,.)\|^2+1)ds.
\end{align}
Then combination (\ref{4.17}) and (\ref{4.20}), we deduce that for all $t\geq T$,
\begin{align} \label{}
\int_{\tau-1}^\tau &e^{\lambda(s-\tau)}z^{2-p}(s,\vartheta_{-\tau}\omega)\|v(s,\tau-t,\vartheta_{-\tau}\omega, v_0)\|_{2p-2}^{2p-2}ds\notag\\
&\leq c z^{-p}(-\tau,\omega)\Big(b(\omega)\int_{-\infty}^{0}e^{\lambda s}z^2(s,\omega)(\|g(s+\tau,.)\|^2+1)ds\notag\\
&+\int_{-\infty}^{0} e^{\lambda s} z^p(s, \omega)(\|g(s+\tau,.)\|^2+1)ds\Big)\notag
\end{align} from which and  (\ref{EE}) it follows that for all $t\geq T$,
\begin{align} \label{4.211}
e^{-\lambda}\int_{\tau-1}^\tau e^{\lambda(s-\tau)} &z^{4-2p}(s,\vartheta_{-\tau}\omega)\|v(s,\tau-t,\vartheta_{-\tau}\omega, v_0)\|^{2p-2}_{2p-2}ds\notag\\
&\leq\int_{\tau-1}^\tau e^{2\lambda(s-\tau)}z^{4-2p}(s,\vartheta_{-\tau}\omega)\|v(s,\tau-t,\vartheta_{-\tau}\omega, v_0)\|^{2p-2}_{2p-2}ds\notag\\
&=z^{p-2}(-\tau,\omega)\int_{\tau-1}^\tau e^{2\lambda(s-\tau)} z^{2-p}(s-\tau,\omega)z^{2-p}(s,\vartheta_{-\tau}\omega)\|v(s)\|^{2p-2}_{2p-2}ds\notag\\
&\leq z^{p-2}(-\tau,\omega)F^{2-p}\int_{\tau-1}^\tau e^{\lambda(s-\tau)} z^{2-p}(s,\vartheta_{-\tau}\omega)\|v(s)\|^{2p-2}_{2p-2}ds\notag\\
&\leq cz^{-2}(-\tau,\omega)F^{2-p}\Big(b(\omega)\int_{-\infty}^{0}e^{\lambda s}z^2(s,\omega)(\|g(s+\tau,.)\|^2+1)ds\notag\\
&+\int_{-\infty}^{0} e^{\lambda s} z^p(s, \omega)(\|g(s+\tau,.)\|^2+1)ds\Big).
\end{align}

For the estimate of the derivative $v_t$ in $L^2_{loc}(\mathbb{R},L^2(\mathbb{R}^N))$, we multiply (\ref{pr1}) by $v_t$ and integrate
over $\mathbb{R}^N$ to produce that
\begin{align}\label{}
&\ \ \ \ \ \ \ \|v_t\|^2+\frac{1}{2}\frac{d}{dt}(\lambda\|v\|^2+\|\nabla v\|^2)\notag\\
&\ \ \ \ \ \ =z(t,\omega)\int_{\mathbb{R}^N}f(x,z^{-1}v)v_tdx+z(t,\omega)\int_{\mathbb{R}^N}gv_tdx\notag\\
&\ \ \ \ \ \ \leq\frac{1}{2}\|v_t\|^2+2\alpha_2^2z^{4-2p}(t,\omega)\|v\|^{2p-2}_{2p-2}+2z^2(t,\omega)\|\psi_2\|^2+z^2(t,\omega)\|g(t,.)\|^2,\notag
\end{align}
\emph{i.e.,} we have
\begin{align}\label{4.22}
&\|v_t\|^2+\frac{d}{dt}(\lambda\|v\|^2+\|\nabla v\|^2)+\lambda(\lambda\|v\|^2+\|\nabla v\|^2)
\notag\\&\ \ \  \ \ \ \ \ \ \ \ \ \ \ \leq cz^{4-2p}(t,\omega)\|v\|^{2p-2}_{2p-2}+4z^2(t,\omega)(\|g(t,.)\|^2+\|\psi_2\|^2)&\notag\\
&\ \ \ \  \ \ \ \ \ \ \ \ \ \ \ \ \ \ \ \ \ \ \ \  +\lambda(\lambda\|v\|^2+\|\nabla v\|^2).
\end{align}
Multiplying (\ref{4.22}) by $e^{\lambda (t-\tau)}$ then integrating about $t$ over $[\tau-1,\tau]$, it give us that, together with
$\omega$  replaced by $\vartheta_{-\tau}\omega$,
\begin{align} \label{4.305}
&\int_{\tau-1}^\tau e^{\lambda(s-\tau)}\|v_s(s,\tau-t,\vartheta_{-\tau}\omega,v_0)\|^2ds\notag\\
&\ \ \  \ \leq c\int_{\tau-1}^\tau e^{\lambda(s-\tau)} z^{4-2p}(s,\vartheta_{-\tau}\omega)\|v(s,\tau-t,\vartheta_{-\tau}\omega,v_0)\|^{2p-2}_{2p-2}ds\notag\\
&\ \ \ \ \ \ \ \ \ +c\int_{\tau-1}^\tau e^{\lambda(s-\tau)}\|v(s,\tau-t,\vartheta_{-\tau}\omega,v_0)\|_{H^1}^2ds\notag\\
&\ \ \ \ \ \ \ \ \ \ \  \ \ \ \ +c\int_{\tau-1}^\tau e^{\lambda s} z^2(s,\vartheta_{-\tau})(\|g(s,.)\|^2+1)ds\notag\\
&\ \ \ \  \ \ \ \ \ \ \  \ \ \ \ \ \ \ +c\|v(\tau-1,\tau-t,\vartheta_{-\tau}\omega,v_0)\|_{H^1}^2.
\end{align}
Then by applying Lemma 4.1 and connection with (\ref{4.211}) and (\ref{4.305}), we deduce that there exists $T=T(\tau,\omega,B)\geq2$ such that for all $t\geq T$,
\begin{align} \label{}
\int_{\tau-1}^\tau e^{\lambda(s-\tau)}&\|v_s(s,\tau-t,\vartheta_{-\tau}\omega,v_0)\|^2ds\notag\\
&\leq c(F^{2-p}b(\omega)+1) L_1(\tau,\omega,\varepsilon)
+cz^{-2}(-\tau,\omega)F^{2-p}\int_{-\infty}^{0} e^{\lambda s} z^p(s, \omega)(\|g(s+\tau,.)\|^2+1)ds\Big).\notag
\end{align}
 This completes the proof.
$\ \ \ \ \ \ \ \ \ \ \ \ \ \ \ \ \ \  \Box$
 \\

We now can prove the $H^1$-tail estimate of solutions of problem (\ref{pr1})-(\ref{pr2}), which is one crucial condition for proving the asymptotic compactness
in $H^1(\mathbb{R}^N)$.\\

\textbf{Lemma 4.4.} \emph{Assume that (\ref{a1})-(\ref{a5}) hold. Given  $\tau\in\mathbb{R}, \omega\in\Omega$ and  $B=\{B(\tau,\omega);\tau\in\mathbb{R},\omega\in\Omega\}\in\mathcal{D}$, then for every $\epsilon>0$, there exist two constants $T=T(\tau,\omega, \epsilon,B)\geq2$ and $R=R(\tau,\omega,\epsilon)>1$
such that the weak solution $v$ of problem (\ref{pr1})-(\ref{pr2}) satisfies that for all $t\geq T$,
\begin{align}\label{}
&\int\limits_{|x|\geq R}\Big(|v(\tau,\tau-t,\vartheta_{-\tau}\omega, z(\tau-t,\vartheta_{-\tau}\omega)u_0)|^2 +|\nabla v(\tau,\tau-t,\vartheta_{-\tau}\omega, z(\tau-t,\vartheta_{-\tau}\omega)u_0)|^2\Big)dx\leq\epsilon,\notag
\end{align} \\
where $u_0\in B(\tau-t,\vartheta_{-t}\omega)$ and $R,T$ are independent of $\varepsilon$.}\\

\emph{Proof}\ \ Given $\xi$ being defined in (4.15), we multiply (\ref{pr1})
by $-\xi\Delta v$ and integrate  over $\mathbb{R}^N$ to find that
\begin{align}\label{4.25}
\frac{1}{2}\frac{d}{dt}\int\limits_{\mathbb{R}^N}\xi|\nabla v|^2dx&+\int\limits_{\mathbb{R}^N}(\nabla\xi.\nabla v)v_tdx\notag\\&+\lambda\int\limits_{\mathbb{R}^N}\xi|\nabla v|^2dx+\lambda\int\limits_{\mathbb{R}^N}(\nabla\xi.\nabla v)vdx+\int\limits_{\mathbb{R}^N}\xi|\Delta v|^2dx\notag\\
&=-z(t,\omega)\int\limits_{\mathbb{R}^N}f(x,z^{-1}v)\xi\Delta vdx-z(t,\omega)\int\limits_{\mathbb{R}^N}g\xi\Delta vdx.
\end{align}
Now, we estimate each term in (\ref{4.25}) as follows. First  we have
\begin{align}\label{4.26}
\Big|\int\limits_{\mathbb{R}^N}(\nabla\xi.\nabla v)v_tdx&+\lambda\int\limits_{\mathbb{R}^N}(\nabla\xi.\nabla v)vdx\Big|\notag\\&=\Big|\int\limits_{\mathbb{R}^N}(v_t+\lambda v)(\frac{2x}{k^2}.\nabla v)\xi^\prime dx\Big|\leq \frac{c}{k}(\|v_t\|^2+\|v\|_{H^1}^2),
\end{align} where and in the following the constant $c$ is independent of $k$ and $\varepsilon$.
For the  nonlinearity in (\ref{4.25}),  we see that
\begin{align}\label{4.27}
-z\int\limits_{\mathbb{R}^N}f(x,z^{-1}v)\xi \Delta vdx&=z\int\limits_{\mathbb{R}^N}f(x,z^{-1}v)(\nabla\xi.\nabla v)dx+z\int\limits_{\mathbb{R}^N}(\frac{\partial}{\partial x} f(x,z^{-1}v).\nabla v)\xi dx\notag\\
&+\int\limits_{\mathbb{R}^N}\frac{\partial}{\partial u}f(x,z^{-1}v)|\nabla v|^2\xi dx.
\end{align}
On the other hand, by using\ (\ref{a2}), (\ref{a3}) and  (\ref{a4}), respectively, we calculate that
\begin{align}\label{4.28}
\Big|z\int\limits_{\mathbb{R}^N}f(x,z^{-1}v)(\nabla\xi.\nabla v)dx\Big|&\leq \frac{z\sqrt{2}C_1}{k}\int\limits_{k\leq |x|\leq \sqrt{2}k}|f(x,z^{-1}v)||\nabla v|dx\notag\\
&\leq \frac{c}{k}(z^{4-2p}\|v\|_{2p-2}^{2p-2}+z^2\|\psi_2\|^2+\|\nabla v\|^2) ,
\end{align}
\begin{align}\label{4.29}
\int\limits_{\mathbb{R}^N}\frac{\partial}{\partial u}f(x,z^{-1}v)|\nabla v|^2\xi dx\leq \alpha_3\int\limits_{\mathbb{R}^N}\xi|\nabla v|^2 dx,
\end{align}
and
\begin{align}\label{4.30}
\Big|z\int\limits_{\mathbb{R}^N}(\frac{\partial}{\partial x} f(x,z^{-1}v).\nabla v)\xi dx\Big|
&\leq \Big|z\int\limits_{\mathbb{R}^N}|\psi_3||\nabla v|\xi dx\Big|\notag\\&\leq \frac{\lambda}{2}\int\limits_{\mathbb{R}^N}\xi|\nabla v|^2 dx+cz^2\int\limits_{\mathbb{R}^N}\xi|\psi_3|^2dx.
\end{align}
Then it follows from \ref{4.27})-(\ref{4.30}) that
\begin{align}\label{4.31}
-z\int\limits_{\mathbb{R}^N}f(x,z^{-1}v)\xi \Delta vdx&\leq \frac{c}{k}(z^{4-2p}\|v\|_{2p-2}^{2p-2}+z^2\|\psi_2\|^2+\|\nabla v\|^2)\notag\\
&+\frac{\lambda}{2}\int\limits_{\mathbb{R}^N}\xi|\nabla v|^2 dx+cz^2\int\limits_{\mathbb{R}^N}\xi|\psi_3|^2dx+\alpha_3\int\limits_{\mathbb{R}^N}\xi|\nabla v|^2 dx.
\end{align}
For the last term on the right hand side of (\ref{4.25}),  we have
\begin{align}\label{4.32}
\Big|z\int\limits_{\mathbb{R}^N}g\xi \Delta vdx\Big|&\leq \frac{1}{2}\int\limits_{\mathbb{R}^N}\xi |\Delta v|^2dx+\frac{1}{2\lambda}z^2\int\limits_{\mathbb{R}^N}\xi |{g}|^2dx.
\end{align}
Then  we incorporate (\ref{4.26}) and (\ref{4.31})-(\ref{4.32}) into (\ref{4.25}) to find that
\begin{align}\label{4.33}
\frac{d}{dt}\int\limits_{\mathbb{R}^N}\xi |\nabla v|^2dx+\lambda\int\limits_{\mathbb{R}^N}\xi |\nabla v|^2dx&\leq \frac{c}{k}(\|v_t\|^2+\|v\|_{H^1}^2+z^{4-2p}\|v\|_{2p-2}^{2p-2}+z^2\|\psi_2\|^2)\notag\\
&+2\alpha_3\int\limits_{\mathbb{R}^N}\xi|\nabla v|^2 dx+cz^2\int\limits_{\mathbb{R}^N}\xi(|\psi_3|^2+|g|^2)dx.
\end{align}
 Applying Lemma 5.1 in \cite{Zhao0} to (\ref{4.33}) over $[\tau-1,\tau]$, we  find that, along with $\omega$ replaced by $\vartheta_{-\tau}\omega$,
\begin{align}\label{4.34}
&\ \ \ \ \ \ \ \ \ \ \int\limits_{\mathbb{R}^N}\xi |\nabla v(\tau,\tau-t,\vartheta_{-\tau}\omega,v_0)|^2dx\notag\\
&\leq \frac{c}{k}\int_{\tau-1}^\tau e^{\lambda(s-\tau)}(\|v_s(s)\|^2+\|v(s)\|_{H^1}^2+z^{4-2p}(s,\vartheta_{-\tau}\omega)\|v(s)\|_{2p-2}^{2p-2}\notag\\
&\  \ \ \ \ \ \ \ \ \ +z^2(s,\vartheta_{-\tau}\omega)\|\psi_2\|^2)ds+c\int_{\tau-1/2}^\tau e^{\lambda(s-\tau)}\int\limits_{|x|\geq k}|\nabla v(s)|^2 dx\notag\\
& \ \ \ \ \ \ \ \ \ \ \ \ \ \ \ +cz^{-2}(\tau,\omega)\int_{-\infty}^0 e^{\lambda s} z^2(s,\omega)\int\limits_{|x|\geq k}(|\psi_3|^2+|g(s+\tau,x)|^2)dxds,
\end{align} where $v(s)=v(s,\tau-t,\vartheta_{-\tau}\omega,z(\tau-t,\vartheta_{-\tau}\omega)u_0)$.
Our task in the following is to show that
each term on the right hand side of (\ref{4.34}) vanishes. First,
by Lemma 4.2,  there are two constants $T_1=T_1(\tau,\omega,B,\epsilon)\geq2$ and $R_1=R_1(\tau,\omega,\epsilon)>1$ such that
for all $t\geq T_1$ and $k\geq R_1$,
\begin{align} \label{4.36}
c\int_{\tau-1}^\tau e^{\lambda(s-\tau)}\int\limits_{|x|\geq k}|\nabla v(s,\tau-t,\vartheta_{-\tau}\omega,v_0)|^2 dxds\leq \frac{\epsilon}{6}.
\end{align}
By Lemma 4.1, it follows that there exist $T_2=T_2(\tau,\omega,B)\geq1$ and $R_2=R_2(\tau,\omega,\epsilon)\geq2$ such that  for all $t\geq T_2$ and $k\geq R_2$,
\begin{align} \label{4.38}
\frac{c}{k}\int_{\tau-1}^\tau e^{\lambda(s-\tau)}\|v(s,\tau-t,\vartheta_{-\tau}\omega,v_0)\|_{H^1}^2ds\leq \frac{\epsilon}{6}.
\end{align}
By Lemma 4.3,  there exist $T_3=T_3(\tau,\omega,B)\geq2$ and $R_3=R_3(\tau,\omega,\epsilon)>1$ such that  for all $t\geq T_3$ and  $k\geq R_3$,
\begin{align} \label{4.39}
&\frac{c}{k}\int_{\tau-1}^\tau e^{\lambda(s-\tau)}z^{4-2p}(s,\vartheta_{-\tau}\omega)\|v(s,\tau-t,\vartheta_{-\tau}\omega,v_0)\|_{2p-2}^{2p-2}ds
\leq \frac{c}{k}L_2(\tau,\omega,\varepsilon)\leq \frac{\epsilon}{6},
\end{align}
and
\begin{align} \label{4.40}
\frac{c}{k}\int_{\tau-1}^\tau e^{\lambda(s-\tau)}\|v_s(s,\tau-t,\vartheta_{-\tau}\omega,v_0)\|^2ds\leq \frac{c}{k}L_3(\tau,\omega,\varepsilon)\leq \frac{\epsilon}{6}.
\end{align}
Similar to (\ref{4.11}), we deduce that there exist $R_4=R_4(\tau,\omega,\epsilon)$ such that for all $k\geq R_4$,
\begin{align} \label{4.37}
cz^{-2}(\tau,\omega)\int_{-\infty}^0e^{\lambda s} z^2(s,\omega)\int\limits_{|x|\geq k}(|\psi_3|^2+|g(s+\tau,x)|^2)dxds\leq \frac{\epsilon}{6}.
\end{align}
Obviously,  there exists $R_5=R_5(\tau,\omega,\epsilon)$ such that t for all  $k\geq R_5$,
\begin{align} \label{4.41}
&\frac{c}{k}\int_{\tau-1}^\tau  e^{\lambda(s-\tau)}z^2(s,\vartheta_{-\tau}\omega)\|\psi_2\|^2ds
\leq\frac{c}{k}\|\psi_2\|^2z^{-2}(-\tau,\omega)\int_{\infty}^0  e^{\lambda s}z^2(s, \omega)ds\leq \frac{\epsilon}{6},
\end{align} where $\int_{\infty}^0  e^{\lambda s}z^2(s, \omega)ds<+\infty$. Finally,
take
$$T=\{T_1,T_2,T_3\},\ \ \   {R}=\max\{R_1,R_2,R_3,R_4,R_5\}.$$
 It is obvious that $R$ and $T$ are independent of the intension $\varepsilon$. Then we combine (\ref{4.36})-(\ref{4.41}) into (\ref{4.34}) to get that for all $t\geq T$ and $k\geq R$,
 \begin{align}  \label{}
\int\limits_{|x|\geq \sqrt{2} k} |\nabla v(\tau,\tau-t,\vartheta_{-\tau}\omega,v_0)|^2dx\leq \epsilon.\notag
\end{align}
  Then connection with  Lemma  4.2, the desired result is achieved.
$\ \ \ \ \ \ \ \ \ \ \ \ \ \ \ \ \ \  \Box$\\

\subsection{Estimate of the truncation of solutions in $L^{2p-2}$ }

Given $u$ the solution of problem (\ref{eq1})-(\ref{eq2}), for each fixed $\tau\in\mathbb{R},\omega\in\Omega$,
we write $M=M(\tau,\omega)>1$ and
\begin{align} \label{}
\mathbb{R}^N(|u(\tau,\tau-t,\vartheta_{-\tau}\omega,u_0)|\geq M)=\{x\in \mathbb{R}^N; |u(\tau,\tau-t,\vartheta_{-\tau}\omega,u_0)|\geq M|\}.\notag
\end{align}

We introduce the trunctation version of solutions of problem (\ref{pr1})-(\ref{pr2}). Let  $(v-M)_+$ be the positive part of $v-M$,\emph{ i.e.},
$$ (v-M)_+=\left\{
       \begin{array}{ll}
    v-M, \ \ \mbox{if}\ v> M;\\
  0,\ \ \ \ \ \ \ \  \ \mbox{if}\ v\leq M.
       \end{array}
      \right.$$

 The next lemma   show that the absolute value  $|u|$ vanishes  in $L^{2p-2}$-norm on the state domain
$\mathbb{R}^N(|u(\tau,\tau-t,\vartheta_{-\tau}\omega),u_0)|\geq M)$ for $M$ large enough, which is the second crucial condition
for proving the asymptotic compactness of solutions in $H^1(\mathbb{R}^N)$.\\

\textbf{Lemma 4.5.}\emph{ Assume that (\ref{a1})-(\ref{a5}) hold. Given  $\tau\in\mathbb{R}, \omega\in\Omega$ and  $B=\{B(\tau,\omega);\tau\in\mathbb{R},\omega\in\Omega\}\in\mathcal{D}$, then for any $\eta>0$, there exist  constants $M=M(\tau,\omega,\eta,B)>1$  and $T=T(\tau,\omega,B)\geq 2$
such that the solution $u_\varepsilon$ of problem (\ref{eq1})-(\ref{eq2}) satisfies that for all $t\geq T$ and  all $\varepsilon\in(0,1]$ and $u_0\in B(\tau-t,\vartheta_{-t}\omega)$,
\begin{align} \label{}
\int_{\tau-1}^\tau e^{\varrho(s-\tau)} \int_{\mathbb{R}^N(|v(s,\tau-t,\vartheta_{-\tau}\omega,z(\tau-t,\vartheta_{-\tau}\omega)u_0)|\geq M)}|v(s,\tau-t,\vartheta_{-\tau}\omega,z(\tau-t,\vartheta_{-\tau}\omega)u_0)|^{2p-2} dxds\leq \eta,\notag
\end{align} where $p>2$ and $M,T$ are independent of  $\varepsilon$  and
$$\varrho=\varrho(\tau,\omega,M)=\alpha_1 E^{2-p}e^{-(p-2)|\omega(-\tau)|}M^{p-2}.$$}

\emph{Proof} \ \
First,  we replace  $\omega$ by $\vartheta_{-\tau}\omega$ in (\ref{pr1}) and see that
$$
v=v(s)=:v(s,\tau-t,\vartheta_{-\tau}\omega,v_0),\ \  s\in[\tau-1,\tau],
$$ is a solution of the following SPDE,
\begin{align} \label{p011}
\frac{dv}{ds}+\lambda v-\Delta v=\frac{z(s-\tau,\omega)}{z(-\tau,\omega)}f(x,u)+\frac{z(s-\tau,\omega)}{z(-\tau,\omega)}{g}(s,x),
\end{align} with the initial data $v_0=z(\tau-t,\vartheta_{-\tau}\omega)u_0$ and $u_0\in B(\tau-t,\vartheta_{-t}\omega)$.

We  multiply (\ref{p011}) by  $(v-M)_+^{p-1}$  and integrate over $\mathbb{R}^N$ to get that for every $s\in[\tau-1,\tau]$,
\begin{align} \label{p01}
&\frac{1}{p}\frac{d}{ds}\int_{\mathbb{R}^N} (v-M)_+^{p}dx+\lambda\int_{\mathbb{R}^N} v(v-M)_+^{p-1}dx-\int_{\mathbb{R}^N}\Delta v(v-M)_+^{p-1}dx
\notag\\&=\frac{z(s-\tau,\omega)}{z(-\tau,\omega)} \int_{\mathbb{R}^N}f(x,u)(v-M)_+^{p-1}dx+\frac{z(s-\tau,\omega)}{z(-\tau,\omega)} \int_{\mathbb{R}^N}{g}(s,x)(v-M)_+^{p-1}dx.
\end{align}
We now have to estimate every term in (\ref{p01}). First, it is obvious that
\begin{align} \label{p02}
-\int_{\mathbb{R}^N}\Delta v(v-M)_+^{p-1}dx=(p-1)\int_{\mathbb{R}^N}(v-M)_+^{p-2}|\nabla
v|^2dx\geq0,
\end{align}
\begin{align} \label{p03}
\lambda\int_{\mathbb{R}^N} v(v-M)_+^{p-1}dx\geq \lambda\int_{\mathbb{R}^N} (v-M)_+^{p}dx.
\end{align}
If $v> M$, then $u=z^{-1}(s,\vartheta_{-\tau}\omega)v>0$, and thus
by (\ref{a1}) and (\ref{EE}), we find that for every $s\in[\tau-1,\tau]$,
\begin{align} \label{}
&\ \  \ \ \ \ \ \ \ f(x,u)\notag\\
&\leq -\alpha_1\Big(\frac{z(s-\tau,\omega)}{z(-\tau,\omega)}\Big)^{1-p}|v|^{p-1}+\frac{z(s-\tau,\omega)}{z(-\tau,\omega)}\frac{\psi_1(x)}{v}\notag\\
&\leq -\frac{1}{2}\alpha_1\Big(\frac{E}{z(-\tau,\omega)}\Big)^{1-p}M^{p-2}(v-M)-\frac{1}{2}\alpha_1\Big(\frac{E}{z(-\tau,\omega)}\Big)^{1-p}(v-M)^{p-1}
+\frac{F}{z(-\tau,\omega)}|\psi_1(x)|(v-M)^{-1},\notag
\end{align}
 by which we find that
\begin{align} \label{p04}
&\ \ \  \ \ \ \ \ \ \ \ \ \ \ \ \ \ \  \frac{z(s-\tau,\omega)}{z(-\tau,\omega)} \int_{\mathbb{R}^N}f(x,u)(v-M)_+^{p-1}dx\notag\\&
\leq-\frac{1}{2}\alpha_1\Big(\frac{E}{z(-\tau,\omega)}\Big)^{2-p}M^{p-2}\int_{\mathbb{R}^N}(v-M)_+^{p}dx-
\frac{1}{2}\alpha_1\Big(\frac{E}{z(-\tau,\omega)}\Big)^{2-p}\int_{\mathbb{R}^N}(v-M)_+^{2p-2}dx\notag\\
&\ \  \ \ \ \ \ \ \ \ \ \ \ \ \ \ \ \ \ +\Big(\frac{F}{z(-\tau,\omega)}\Big)^{2}\int_{\mathbb{R}^N}|\psi_1(x)|(v-M)_+^{p-2}dx\notag\\
&\leq-\frac{1}{2}\alpha_1\Big(\frac{E}{z(-\tau,\omega)}\Big)^{2-p}M^{p-2}\int_{\mathbb{R}^N}(v-M)_+^{p}dx-\frac{1}{2}\alpha_1\Big(\frac{E}{z(-\tau,\omega)}
\Big)^{2-p}\int_{\mathbb{R}^N}(v-M)_+^{2p-2}dx\notag\\
& \ \ \ \ \  \ \ \ \ \ \ \ \  \ \ \ \ +\frac{1}{2}\lambda\int_{\mathbb{R}^N}(v-M)_+^{p}dx+c\Big(\frac{F}{z(-\tau,\omega)}\Big)^{p}\int_{\mathbb{R}^N(v\geq M)}|\psi_1(x)|^{p/2}dx.
\end{align}
The second term on the right hand side of (\ref{p01}) is estimated as
\begin{align} \label{p05}
\frac{F}{z(-\tau,\omega)}\Big| \int_{\mathbb{R}^N}{g}(s,x)(v(s)-M)_+^{p-1}dx\Big|&\leq \frac{1}{4}\alpha_1\Big(\frac{E}{z(-\tau,\omega)}
\Big)^{2-p}\int_{\mathbb{R}^N}(v-M)_+^{2p-2}dx\notag\\
&+\frac{1}
{\alpha_1}\Big(\frac{F}{z(-\tau,\omega)}\Big)^2\Big(\frac{E}{z(-\tau,\omega)}
\Big)^{p-2}\int_{\mathbb{R}^N(v(s)\geq M)}{g}^2(s,x)dx\notag\\
&\leq \frac{1}{4}\alpha_1\Big(\frac{E}{z(-\tau,\omega)}
\Big)^{2-p}\int_{\mathbb{R}^N}(v-M)_+^{2p-2}dx\notag\\
&+\frac{1}
{\alpha_1}\Big(\frac{F}{z(-\tau,\omega)}
\Big)^{p}\int_{\mathbb{R}^N(v(s)\geq M)}{g}^2(s,x)dx.
\end{align}
Combination (\ref{p01})-(\ref{p05}), we obtain that
\begin{align} \label{p060}
\frac{d}{ds}\int_{\mathbb{R}^N} (v(s)-M)_+^{p}dx&+
\alpha_1\Big(\frac{E}{z(-\tau,\omega)}\Big)^{2-p}M^{p-2}\int_{\mathbb{R}^N}(v(s)-M)_+^{p}dx\notag\\
&+\alpha_1\Big(\frac{E}{z(-\tau,\omega)}
\Big)^{2-p}\int_{\mathbb{R}^N}(v-M)_+^{2p-2}dx\notag\\
&\leq c\Big(\frac{F}{z(-\tau,\omega)}\Big)^{p}\Big(\|{g}(s,.)\|^2+
\|\psi_1\|_{p/2}^{p/2}\Big),
\end{align}
where the positive constant $c$  is independent of $\varepsilon,\tau,\omega$ and $M$.  Note that for each $\tau\in \mathbb{R}$ and $\varepsilon\in(0,1]$,
\begin{align} \label{TT}
 e^{-|\omega(-\tau)|} \leq z^{-1}(-\tau,\omega)=e^{\varepsilon \omega(-\tau)}\leq e^{|\omega(-\tau)|}.
 \end{align}
For convenience, we put
$$
 \varrho=\varrho(\tau,\omega,M)=\alpha_1 E^{2-p}e^{-(p-2)|\omega(-\tau)|}M^{p-2},\ \ d=d(\tau,\omega)=\alpha_1E^{2-p}e^{-(p-2)|\omega(-\tau)|}
$$ Then (\ref{p060}) is rewrote as
\begin{align} \label{p06}
\frac{d}{ds}\int_{\mathbb{R}^N} (v(s)-M)_+^{p}dx&+
\varrho\int_{\mathbb{R}^N}(v(s)-M)_+^{p}dx+d\int_{\mathbb{R}^N}(v-M)_+^{2p-2}dx\notag\\
&\leq cF^pe^{p|\omega(-\tau)|}\Big(\|{g}(s,.)\|^2+1\Big),
\end{align} where $s\in[\tau-1,\tau]$ and $\varrho,E,F$ are independent of $\varepsilon$.
By using Lemma 5.1 in \cite{Zhao0}  to (\ref{p06}) over $[\tau-1,\tau]$, we
find that
 \begin{align} \label{p07}
d\int_{\tau-1}^\tau\int_{\mathbb{R}^N}e^{\varrho(s-\tau)}(v(s)-M)_+^{2p-2}dxds
&\leq \int_{\tau-1}^\tau e^{\varrho(s-\tau)}\int_{\mathbb{R}^N} \Big(v(s,\tau-t,\vartheta_{-\tau}\omega,v_0)-M\Big)_+^{p}dxds\notag\\
&+ cF^pe^{p|\omega(-\tau)|}\int_{\tau-1}^\tau e^{\varrho(s-\tau)} \Big(\|g(s,.)\|^2+1\Big)ds.
\end{align}
First by (\ref{4.17}), there exists $T=T(\tau,\omega,B)\geq2$ such that for all $t\geq T$,
\begin{align} \label{p077}
\int_{\tau-1}^\tau e^{\varrho(s-\tau)}\int_{\mathbb{R}^N} \Big(v(s,\tau-t,\vartheta_{-\tau}\omega,v_0)-M\Big)_+^{p}dxds\leq N(\tau,\omega,\varepsilon)\frac{1}{\varrho}\rightarrow0,
\end{align} as $\varrho\rightarrow+\infty$, where $N(\tau,\omega,\varepsilon)$ is defined by the right hand side of (\ref{4.17}). We then need to show  the second term on the right hand side of (\ref{p07}) is also small as $\varrho\rightarrow+\infty$.  Indeed,
choosing $\varrho>\delta$ and taking $\varsigma\in(0,1)$, we have
\begin{align} \label{}
\int_{\tau-1}^\tau e^{\varrho(s-\tau)} \Big(\|g(s,.)\|^2+1\Big)ds&=\int_{\tau-1}^{\tau-\varsigma} e^{\varrho(s-\tau)} (\|g(s,.)\|^2+1)ds+\int_{\tau-\varsigma}^{\tau} e^{\varrho(s-\tau)}(\|g(s,.)\|^2+1)ds\notag\\
&=e^{-\varrho\tau}\int_{\tau-1}^{\tau-\varsigma} e^{(\varrho-\delta)s}e^{\delta s}(\|g(s,.)\|^2+1)ds+e^{-\varrho\tau}\int_{\tau-\varsigma}^{\tau} e^{\varrho s}(\|g(s,.)\|^2+1)ds\notag\\
&\leq e^{-\varrho\varsigma} e^{\delta(\varsigma-\tau)}\int_{-\infty}^{\tau} e^{\delta s}(\|g(s,.)\|^2+)ds+\int_{\tau-\varsigma}^{\tau}(\|g(s,.)\|^2+1)ds.\notag
\end{align}
By (\ref{a5}), the first term above vanishes as $\varrho\rightarrow+\infty$, and by $g\in L^2_{loc}(\mathbb{R}, L^2(\mathbb{R}^N))$  we can choose $\varsigma$ small enough such that the second term is small. Then when $\varrho\rightarrow +\infty$, we have
  \begin{align} \label{p0771}
 cF^pe^{p|\omega(-\tau)|}\int_{\tau-1}^\tau e^{\varrho(s-\tau)} \Big(\|g(s,.)+\|^2+1\Big)ds\rightarrow0.
\end{align} Since if $M\rightarrow+\infty$, then $\varrho\rightarrow +\infty$,  so by (\ref{p07})-(\ref{p0771}), we know that for  $M\rightarrow+\infty$,
 \begin{align} \label{p0772}
 \int_{\tau-1}^\tau e^{\varrho(s-\tau)}\int_{\mathbb{R}^N}(v(s)-M)_+^{2p-2}dxds\rightarrow 0.
 \end{align}
 Note that $v-M\geq \frac{v}{2}$ for $v\geq2M$. Then by (\ref{p0772}) it gives that
 \begin{align} \label{}
 \int_{\tau-1}^\tau e^{\varrho(s-\tau)}\int_{\mathbb{R}^N(v(s)\geq -2M)}|v(s)|^{2p-2}dxds\rightarrow 0, \notag
 \end{align}as $M\rightarrow+\infty$.
 By a similar argument, we can show that there exists $T=T(\tau,\omega,B)\geq2$ such that for all $t\geq T$,
 \begin{align} \label{}
 \int_{\tau-1}^\tau e^{\varrho(s-\tau)}\int_{\mathbb{R}^N(v(s)\leq 2M)}|v(s)|^{2p-2}dxds\rightarrow 0,\notag
 \end{align} as $M\rightarrow+\infty$. Then we finish the total proof.
 $\ \ \  \ \ \ \ \ \ \Box$\\

\subsection{Asymptotic compactness  on bounded domains}

In this subsection,  by using Lemma 4.5, we prove the asymptotic compactness of the cocyle $\varphi$ defined by (\ref{eq0}) in $H_0^1(\mathcal{O}_R)$ for any $R>0$, where $\mathcal{O}_R=\{x\in \mathbb{R}^N; |x|\leq R\}$. For this purpose, we define $\phi(.)=1-\xi(.)$, where $\xi$ is the cut-off function as in (4.15). Then we know that $0\leq\phi(s)\leq1$, and $\phi(s)=1$ if $s\in[0,1]$ and $\phi(s)=0$ if $s\geq2$.
Fix a positive constant $k$, we define
\begin{align} \label{vuu}
\tilde{v}(t,\tau,\omega,v_0)=\phi(\frac{x^2}{k^2})v(t,\tau,\omega,v_0),\ \  \tilde{u}(t,\tau,\omega,u_0)=\phi(\frac{x^2}{k^2})u(t,\tau,\omega,u_0),
\end{align}
where $v$ is the solution of problem (\ref{pr1})-(\ref{pr2}) and $u$ is the solution of problem ({\ref{eq1}})-({\ref{eq2}}) with $v=z(t,\omega)u$. Then we have
\begin{align} \label{vu}
\tilde{u}(t,\tau,\omega,u_0)=z^{-1}(t,\omega)\tilde{v}(t,\tau,\omega,v_0).
\end{align}
It is obvious that $\tilde{v}$
solves the following equations:
\begin{equation}
\begin{cases} \tilde{v}_t+\lambda\tilde{v}-\Delta \tilde{v}=\phi zf(x,z^{-1}v)+\phi zg-v\Delta \phi-2\nabla \phi.\nabla v,\\
\tilde{v}|_{\partial\mathcal{O}_{k\sqrt{2}}}=0,\\
\tilde{v}(\tau,x)=\tilde{v}_0(x)=\phi v_0(x),
\end{cases}
\end{equation} where $\phi=\phi(\frac{x^2}{k^2})$.

It is well-known that the eigenvalue problem on bounded domains $\mathcal{O}_{k\sqrt{2}}$ with Dirichlet boundary condition:
\begin{equation}
\begin{cases} -\Delta \tilde{v}=\lambda \tilde{v},\\
\tilde{v}|_{\partial\mathcal{O}_{k\sqrt{2}}}=0\notag
\end{cases}
\end{equation}
has a family of orthogonal eigenfunctions $\{e_j\}_{=1}^{+\infty}$ in both $L^2(\mathcal{O}_{k\sqrt{2}})$ and
$H_0^1(\mathcal{O}_{k\sqrt{2}})$ such that the corresponding eigenvalue $\{\lambda_j\}_{j=1}^{+\infty}$ is non-decreasing in $j$.

Let $H_m=\mbox{Span}\{e_1,e_2,...,e_m\}\subset H_0^1(\mathcal{O}_{k\sqrt{2}})$
and $P_m:  H^1_0(\mathcal{O}_{k\sqrt{2}})\rightarrow H_m$ be the canonical
projector and $I$ be the identity. Then for every $\tilde{u}\in
H^1_0(\mathcal{O}_{k\sqrt{2}})$, $\tilde{u}$ has a unique decomposition: $\tilde{u}=\tilde{u}_1+\tilde{u}_2$,
where $\tilde{u}_1=P_m\tilde{u}\in H_m$
 and $\tilde{u}_2=(I-P_m)\tilde{u}\in H_m^{\bot}$, \emph{i.e.}, $ H^1_0(\mathcal{O}_{k\sqrt{2}})=H_m\oplus H_m^{\bot}$.\\

\textbf{Lemma 4.6.} \emph{ Assume that (\ref{a1})-(\ref{a5}) hold. Given  $\tau\in\mathbb{R}, \omega\in\Omega$ and  $B=\{B(\tau,\omega);\tau\in\mathbb{R},\omega\in\Omega\}\in\mathcal{D}$, then for every $\epsilon>0$, there
are $N_0=N_0(\tau,\omega, k,\epsilon)\in Z^+$ and
$T=T(\tau,\omega,B,\epsilon)\geq 2$ such that for all $t\geq T$ and
$m> N_0$,
\begin{align} \label{}
 \| (I-P_m)\tilde{u}(\tau,\tau-t,\vartheta_{-\tau}\omega, \tilde{u}_0)\|_{H^1_0(\mathcal{O}_{k\sqrt{2}})}\leq\epsilon,\notag
\end{align}  where $\tilde{u}_0=\phi{u}_0$ with ${u}_0\in B(\tau-t,\vartheta_{-\tau}\omega)$.  Here $\tilde{u}$ is as in (\ref{vu}) and $N,T$ are independent of $\varepsilon$.}
\\

\emph{Proof}\ \  By (\ref{vu}), we start at the estimate of $\tilde{v}$. For $\tilde{v}\in H^1_0(\mathcal{O}_{k\sqrt{2}})$,
we write $\tilde{v}=\tilde{v}_1+\tilde{v}_2$ where $\tilde{v}_1=P_m\tilde{v}$ and
$\tilde{v}_2=(I-P_m)\tilde{v}$.
 Then naturally, we have a splitting about $\tilde{u}=\tilde{u}_1+\tilde{u}_2$  where $\tilde{u}_1=P_m\tilde{u}$ and
$\tilde{u}_2=(I-P_m)\tilde{u}$.
 Multiplying (4.47) by $\Delta \tilde{v}_2$ we get that
\begin{align} \label{4.44}
&\frac{1}{2}\frac{d}{dt}\|\nabla \tilde{v}_2\|_{L^2(\mathcal{O}_{k\sqrt{2}})}^2+\lambda\|\nabla \tilde{v}_2\|_{L^2(\mathcal{O}_{k\sqrt{2}})}^2+\|\Delta\tilde{v}_2\|_{L^2(\mathcal{O}_{k\sqrt{2}})}^2\notag\\
&=-z\int\limits_{\mathcal{O}_{k\sqrt{2}}}\phi f(x,z^{-1}v)\Delta \tilde{v}_2dx+\int\limits_{\mathcal{O}_{k\sqrt{2}}}(\phi zg-v\Delta \phi-2\nabla \phi.\nabla v)\Delta \tilde{v}_2dx.
\end{align}
By (\ref{a2}), we deduce that
\begin{align} \label{4.45}
z\int\limits_{\mathcal{O}_{k\sqrt{2}}}\phi f(x,z^{-1}v)\Delta \tilde{v}_2dx&\leq \frac{1}{4}\|\Delta\tilde{v}_2\|_{L^2(\mathcal{O}_{k\sqrt{2}})}^2+cz^{4-2p}\|v\|_{L^{2p-2}(\mathcal{O}_{k\sqrt{2}})}^{2p-2}+z^2\|\psi_2\|^2.
\end{align}
On the other hand,
\begin{align} \label{4.46}
\int\limits_{\mathcal{O}_{k\sqrt{2}}}(\phi zg-v\Delta \phi-2\nabla \phi.\nabla v)\Delta \tilde{v}_2dx\leq \frac{1}{4}\|\Delta\tilde{v}_2\|_{L^2(\mathcal{O}_{k\sqrt{2}})}^2
+c(z^2\|g\|^2+\|{v}\|^2+\|\nabla{v}\|^2).
\end{align}
Then by (\ref{4.44})-(\ref{4.46}) we find that
\begin{align} \label{}
\frac{d}{dt}\|\nabla \tilde{v}_2\|_{L^2(\mathcal{O}_{k\sqrt{2}})}^2+\|\Delta\tilde{v}_2\|_{L^2(\mathcal{O}_{k\sqrt{2}})}^2
\leq c(z^{4-2p}\|v\|_{L^{2p-2}(\mathcal{O}_{k\sqrt{2}})}^{2p-2}+z^2\|\psi_2\|^2+z^2\|g\|^2+\|{v}\|_{H^1}^2).\notag
\end{align}
from which and connection with the Poincar\'{e}'s inequality $$\|\Delta\tilde{v}_2\|_{L^2(\mathcal{O}_{k\sqrt{2}})}^2\geq \lambda_{m+1}\|\nabla\tilde{v}_2\|_{L^2(\mathcal{O}_{k\sqrt{2}})}^2,$$ it follows that
\begin{align} \label{4.47}
\frac{d}{dt}\|\nabla \tilde{v}_2\|_{L^2(\mathcal{O}_{k\sqrt{2}})}^2&+\lambda_{m+1}\|\nabla\tilde{v}_2\|_{L^2(\mathcal{O}_{k\sqrt{2}})}^2
\notag\\&\leq c(z^{4-2p}\|v\|_{L^{2p-2}(\mathcal{O}_{k\sqrt{2}})}^{2p-2}+z^2\|\psi_2\|^2+z^2\|g\|^2+\|{v}\|_{H^1}^2).
\end{align}
 Applying Lemma 5.1 in \cite{Zhao0} to (\ref{4.47}) over the interval $[\tau-1,\tau]$,  we find that, along with $\omega$  replaced by $\vartheta_{-\tau}\omega$,
\begin{align} \label{4.48}
&\ \  \ \ \ \ \  \ \ \ \ \  \ \  \ \|\nabla \tilde{v}_2(\tau,\tau-t,\vartheta_{-\tau}\omega,\tilde{v}_0)\|_{L^2(\mathcal{O}_{k\sqrt{2}})}^2
\notag\\&\leq \int_{\tau-1}^\tau e^{\lambda_{m+1}(s-\tau)}\|\nabla \tilde{v}_2(s,\tau-t,\vartheta_{-\tau}\omega,\tilde{v}_0)\|_{L^2(\mathcal{O}_{k\sqrt{2}})}^2ds\notag\\
&+c\int_{\tau-1}^\tau e^{\lambda_{m+1}(s-\tau)}z^{4-2p}(s,\vartheta_{-\tau}\omega)\|v(s,\tau-t,\vartheta_{-\tau}\omega,\tilde{v}_0)\|_{L^{2p-2}(\mathcal{O}_{k\sqrt{2}})}^{2p-2}ds\notag\\
&+c\int_{\tau-1}^\tau e^{\lambda_{m+1}(s-\tau)}(z^2(s,\vartheta_{-\tau}\omega)\|\psi_2\|^2+z^2(s,\omega)\|g(s,.)\|^2)ds\notag\\
&+c\int_{\tau-1}^\tau e^{\lambda_{m+1}(s-\tau)}\|{v}(s,\tau-t,\vartheta_{-\tau}\omega,\tilde{v}_0)\|_{H^1}^2ds\notag\\
&\ \  \ \ \ \ \ \ \ \ \ \leq+c\int_{\tau-1}^\tau e^{\lambda_{m+1}(s-\tau)}z^{4-2p}(s,\vartheta_{-\tau}\omega)\|v(s,\tau-t,\vartheta_{-\tau}\omega,\tilde{v}_0)\|_{L^{2p-2}(\mathcal{O}_{k\sqrt{2}})}^{2p-2}ds\notag\\
&\ \ \ \ \ \ \ \ \ \ \ \ +\int_{\tau-1}^\tau e^{\lambda_{m+1}(s-\tau)}\|v(s,\tau-t,\vartheta_{-\tau}\omega,\tilde{v}_0)\|_{H^1}^2ds\notag\\
&\ \  \ \ \ \ \ \ \ \ \ \ +c\int_{\tau-1}^\tau e^{\lambda_{m+1}(s-\tau)}z^2(s,\vartheta_{-\tau}\omega)\Big(\|g(s,.)\|^2+1\Big)ds\notag\\
&\ \ \ \ \ \ \ \ \ \ \ \ \ =I_1+I_2+I_3.
\end{align}
We next to show that $I_1,I_2$ and $I_3$ converge to zero as $m$ increases to infinite.
First  we have
\begin{align} \label{4.50}
I_1&=z^{2p-4}(-\tau,\omega)\int_{\tau-1}^\tau e^{\lambda_{m+1}(s-\tau)}z^{4-2p}(s-\tau,\omega)\|v(s,\tau-t,\vartheta_{-\tau}\omega,\tilde{v}_0)\|_{L^{2p-2}(\mathcal{O}_{k\sqrt{2}})}^{2p-2}ds\notag\\
&\leq z^{2p-4}(-\tau,\omega)F^{4-2p}\int_{\tau-1}^\tau e^{\lambda_{m+1}(s-\tau)}\|v(s,\tau-t,\vartheta_{-\tau}\omega,\tilde{v}_0)\|_{L^{2p-2}(\mathcal{O}_{k\sqrt{2}})}^{2p-2}ds\notag\\
&\leq z^{2p-4}(-\tau,\omega)F^{4-2p}\Big(\int_{\tau-1}^\tau e^{\lambda_{m+1}(s-\tau)}\int_{\mathcal{O}_{k\sqrt{2}}(|v(s)|\geq M)}|v(s,\tau-t,\vartheta_{-\tau}\omega,\tilde{v}_0)|^{2p-2}dxds\notag\\
&+\int_{\tau-1}^\tau e^{\lambda_{m+1}(s-\tau)}\int_{\mathcal{O}_{k\sqrt{2}}(|v(s)|\leq M)}|v(s,\tau-t,\vartheta_{-\tau}\omega,\tilde{v}_0)|^{2p-2}dxds\Big).
\end{align}
By Lemma 4.5  there exist $T_1=T_1(\tau,\omega,B)\geq2$, $M=M(\tau,\omega,B)$  such that for all  $t\geq T_1$,
\begin{align} \label{4.501}
\int_{\tau-1}^\tau e^{\varrho(s-\tau)} \int_{\mathcal{O}_{k\sqrt{2}}(|v(s)|\geq M)}|v(s,\tau-t,\vartheta_{-\tau}\omega,\tilde{v}_0)|^{2p-2} dxds\leq \epsilon.
\end{align}
But $\lambda_{m+1}\rightarrow +\infty$, then there exists $N^\prime=N^\prime(\tau,\omega)>0$ such that  for all $m>N^\prime$, $\lambda_{m+1}>\varrho$. Hence
by (\ref{4.501}) it gives us that  for all  $t\geq T_1$ and $m>N^\prime$ there holds
\begin{align} \label{4.5010}
\int_{\tau-1}^\tau e^{\lambda_{m+1}(s-\tau)} \int_{\mathcal{O}_{k\sqrt{2}}(|v(s)|\geq M)}|v(s,\tau-t,\vartheta_{-\tau}\omega,\tilde{v}_0)|^{2p-2} dxds\leq \epsilon.
\end{align}
For the second term on the right hand side of (\ref{4.50}), since  $\mathcal{O}_{k\sqrt{2}}(|v(s)|\leq M)$ is a bounded domain, then there exists  $N^{\prime\prime}=N^{\prime\prime}(\tau,\omega)>0$ such that
for all $m>N^{\prime\prime}$,
\begin{align} \label{4.502}
\int_{\tau-1}^\tau e^{\lambda_{m+1}(s-\tau)}&\int_{\mathcal{O}_{k\sqrt{2}}(|v(s)|\leq M)}|v(s,\tau-t,\vartheta_{-\tau}\omega,\tilde{v}_0)|^{2p-2}dxds\notag\\
&\leq \frac{1}{\lambda_{m+1}}M^{2p-2}.|(\mathcal{O}_{k\sqrt{2}}(|v(s)|\leq M))|\leq\epsilon,
\end{align} where $|(\mathcal{O}_{k\sqrt{2}}(|v(s)|\leq M))|$ is the finite measure of the bounded domain $\mathcal{O}_{k\sqrt{2}}(|v(s)|\leq M)$.
Put $N_1=\max\{N^{\prime},N^{\prime\prime}\}$. It follows from (\ref{4.50})-(\ref{4.502}) that for all $m>N_1$ and $t\geq T_1$,
\begin{align} \label{4.55}
I_1\leq C_1(\tau,\omega)\epsilon.
\end{align}
By Lemma 4.1, there exists $T_2=T_2(\tau,\omega)$ and $N_2=N_2(\tau,\omega)>0$ such that for all  $m>N_2$ and $t\geq T_2$,
\begin{align} \label{4.56}
I_2\leq \frac{L_1(\tau,\omega,\varepsilon)}{\lambda_{m+1}}\leq\epsilon.
\end{align}
By a same technique as (\ref{p0771}), we can show that there exists $N_3=N_3(\tau,\omega)>0$ such that for all  $m>N_3$,
\begin{align} \label{4.551}
I_3=\int_{\tau-1}^\tau e^{\lambda_{m+1}(s-\tau)}z^2(s,\vartheta_{-\tau}\omega)\Big(\|g(s,.)\|^2+1\Big)ds\leq \epsilon.
\end{align}
Let $N_0=\max\{N_1,N_2,N_3\}$ and $T=\max\{T_1,T_2\}$. Then combination (\ref{4.48}) and (\ref{4.55})-(\ref{4.551}), we get that there exists a finite
constant $\mu=\mu(\tau,\omega)>0$
 such that for all  $m>N_0$ and $t\geq T$,
\begin{align} \label{4.126}
\|\nabla \tilde{v}_2(\tau,\tau-t,\vartheta_{-\tau}\omega, \tilde{v}_0)\|_{L^2(\mathcal{O}_{k\sqrt{2}})}\leq C_1(\tau,\omega)\epsilon.
\end{align}  Then by (\ref{eq00}) and  (\ref{4.126}), we have
\begin{align} \label{}
\|\nabla \tilde{u}_2(\tau,\tau-t,\vartheta_{-\tau}\omega, \tilde{u}_0)\|_{L^2(\mathcal{O}_{k\sqrt{2}})}=z(-\tau,\omega)\|\nabla\tilde{v}_2(\tau,\tau-t,\vartheta_{-\tau}\omega, \tilde{v}_0)\|_{L^2(\mathcal{O}_{k\sqrt{2}})}\leq C_2(\tau,\omega)\epsilon,\notag
\end{align}for all  $m>N_0$ and $t\geq T$, which completes the proof.
 $\ \ \ \ \ \ \Box$\\

 \textbf{Lemma 4.7.} \emph{Assume that (\ref{a1})-(\ref{a5}) hold. Given  $\tau\in\mathbb{R}, \omega\in\Omega$, then for every $k>0$, the sequence $\{\tilde{u}(\tau, \tau-t_n,\vartheta_{-\tau}\omega, \phi(\frac{x^2}{k^2})u_{0,n})\}_{n=1}^\infty$ has a convergent subsequence
 in $H_0^1(\mathcal{O}_{k\sqrt{2}})$ whenever $t_n\rightarrow+\infty$ and $u_{0,n}\in B(\tau-t_n,\vartheta_{-t_n}\omega)$.}\\

\emph{Proof}\ \ Given $\epsilon>0$, by Lemma 4.6, there exist $N_0\in \mathbb{Z}^+$ such that  as $t_n\rightarrow+\infty$
 \begin{align} \label{L101}
\|(I-P_{N_0})\tilde{u}(\tau,\tau-t_n,\vartheta_{-\tau}\omega, \phi(\frac{x^2}{k^2})u_{0,n})\|_{H^1(\mathcal{O}_{k\sqrt{2}})}\leq \epsilon.
\end{align}
 By Lemma 4.1, we deduce that  if $t_n$ large enough,
 \begin{align} \label{L102}
\|P_{N_0}\tilde{u}(\tau,\tau-t_n,\vartheta_{-\tau}\omega, \phi(\frac{x^2}{k^2})u_{0,n})\|_{H^1(\mathcal{O}_{k\sqrt{2}})}\leq L_1(\tau,\omega,\varepsilon).
\end{align}
 Note that $H^1(\mathcal{O}_{k\sqrt{2}})=P_{N_0}H^1(\mathcal{O}_{k\sqrt{2}})+(I-P_{N_0})H^1(\mathcal{O}_{k\sqrt{2}})$, but $P_{N_0}H^1(\mathcal{O}_{k\sqrt{2}})$ is a finite
 dimensional space. Then by (\ref{L102}),  if $n,m$ large enough,
  \begin{align} \label{L103}
\|P_{N_0}\tilde{u}(\tau,\tau-t_n,\vartheta_{-\tau}\omega, \phi(\frac{x^2}{k^2})u_{0,n})-P_{N_0}\tilde{u}(\tau,\tau-t_m,\vartheta_{-\tau}\omega, \phi(\frac{x^2}{k^2})u_{0,m}\|_{H^1(\mathcal{O}_{k\sqrt{2}})}\leq \epsilon.
\end{align}
  Then it is easy to finish the proof by means of (\ref{L101}) and (\ref{L103}) and a standard argument.
$\ \ \ \ \ \ \Box$

\subsection{Existence of pullback attractor in $H^{1}(\mathbb{R}^N)$}

In this subsection, we prove the existences of pullback attractors in $H^{1}(\mathbb{R}^N)$ for problem (\ref{eq1})-(\ref{eq2}) for every $\varepsilon\in(0,1]$.\\

\textbf{ Lemma 4.8.} \emph{Assume that (\ref{a1})-(\ref{a5}) hold . Then the cocycle $\varphi$ defined by (\ref{eq0}) is asymptotically compact
in $H^1(\mathbb{R}^N)$, i.e., for every  $\tau\in\mathbb{R}, \omega\in\Omega$,  the sequence $\{\varphi(t, \tau-t_n,\vartheta_{-t}\omega, u_{0,n})\}_{n=1}^\infty$ has a convergent subsequence
 in $H^1(\mathbb{R}^N)$ whenever $t_n\rightarrow+\infty$ and $u_{0,n}\in B=B(\tau-t_n,\vartheta_{-t_n}\omega)$ with $B\in\mathcal{D}$.}\\

\emph{Proof}\ \  Give $R>0$, denote by $\mathcal{O}^c_{R}=\mathbb{R}^N-\mathcal{O}_{R}$, where $\mathcal{O}_R=\{x\in \mathbb{R}^N; |x|\leq R\}$. By  Lemma 4.4,  for
 any $\epsilon>0$, there exist $R=R(\tau,\omega,\epsilon)>0$ and $N_1=N_1(\tau,\omega,B,\epsilon)\in \mathbb{Z}^+$ such that for all $n\geq N_1$,
\begin{align} \label{n08}
\|v(\tau, \tau-t_n,\vartheta_{-\tau}\omega, z(\tau-t_n,\vartheta_{-\tau}\omega)u_{0,n})\|_{H^1(\mathcal{O}^c_{R})}\leq \frac{\epsilon}{8} e^{-|\omega(-\tau)|},
\end{align} for every $u_{0,n}\in B=B(\tau-t_n,\vartheta_{-t_n}\omega)$.
By (\ref{eq00}) and (\ref{n08}), we have
\begin{align} \label{n0800}
\|u(\tau, \tau-t_n,\vartheta_{-\tau}\omega, z(\tau-t_n,\vartheta_{-\tau}\omega)u_{0,n})\|_{H^1(\mathcal{O}^c_{R})}\leq \frac{\epsilon}{8}.
\end{align}
On the other hand, for this radius $R$, by Lemma 4.7, there exists $N_2=N_2(\tau,\omega,B,\epsilon)\geq N_1$ such that for all $m,n\geq N_2$,
\begin{align} \label{n081}
&\|u(\tau, \tau-t_n,\vartheta_{-\tau}\omega, \phi(\frac{x^2}{R^2})\omega)u_{0,n})-u(\tau, \tau-t_m,\vartheta_{-\tau}\omega, \phi(\frac{x^2}{R^2})u_{0,m})\|_{H^1_0(\mathcal{O}_{R\sqrt{2}})}\leq \frac{\epsilon}{8}.
\end{align} Then the desired result follows from (\ref{n0800}) and (\ref{n081}) by a standard argument.
$\ \ \ \ \ \ \Box$\\

Given $\varepsilon\in (0,1]$, by Lemma 4.1, we deduce that the $\mathcal{D}$-pullback absorbing set $K_\varepsilon$ of $\varphi_{\varepsilon}$ in $L^2(\mathbb{R}^N)$ is defined by
\begin{align}\label{KK1}
K_\varepsilon=\{K_\varepsilon(\tau,\omega)=\{u\in L^2(\mathbb{R}^N); \|u\|\leq L_\varepsilon(\tau,\omega)\};\tau\in\mathbb{R},\omega\in
\Omega\},
\end{align}  where
\begin{align}\label{KK2}
L_\varepsilon(\tau,\omega,\varepsilon)=c\Big(\int_{-\infty}^0e^{\lambda s}e^{-2\varepsilon \omega(s)}(\|g(s+\tau,.)\|^2+1)\Big)^{1/2}.
\end{align}

By Lemma 4.8 and  Theorem 2.6,  we immediately have\\

{\bf Theorem 4.9.} \emph{ Assume that (\ref{a1})-(\ref{a5}) hold. Then for every fixed $\varepsilon\in(0,1]$, the cocycle $\varphi_\varepsilon$ defined by (\ref{eq0}) possesses a
unique $\mathcal{D}$-pullback attractor $\mathcal{A}_{\varepsilon,H^1}=\{\mathcal{A}_{\varepsilon,H^1}(\tau,\omega);\tau\in\mathbb{R},\omega\in \Omega\}$ in $H^1(\mathbb{R}^N)$, given by
\begin{align}\label{}
\mathcal{A}_{\varepsilon,H^1}(\tau,\omega)=\bigcap_{s>0}\overline{\bigcup_{t\geq s} \varphi_\varepsilon(t,\tau-t,\vartheta_{-t}\omega, K_\varepsilon(\tau-t,\vartheta_{-t}\omega))}^{H^1(\mathbb{R}^N)}, \ \
\tau\in \mathbb{R},\omega\in \Omega.\notag
 \end{align}
Furthermore,  $\mathcal{A}_{\varepsilon,H^1}$ is consistent with the $\mathcal{D}$-pullback random attractor $\mathcal{A}_\varepsilon$ in  $L^2(\mathbb{R}^N)$, which is defined as in (\ref{L2}).}

\section{Upper semi-continuity  of pullback attractor in $H^{1}(\mathbb{R}^N)$}

From Theorem 4.9, for every $\varepsilon\in(0,1]$, the cocycle $\varphi_\varepsilon$ admits a common $\mathcal{D}$-pullback  attractor $\mathcal{A}_\varepsilon$ in both $L^2(\mathbb{R}^N)$ and $H^1(\mathbb{R}^N)$, where $\mathcal{D}$ is
defined by (3.11).  From this fact we may investigate the upper semi-continuity of $\mathcal{A}_\varepsilon$ in both $L^2(\mathbb{R}^N)$ and $H^1(\mathbb{R}^N)$.
Note that \cite{Wang1} only proved the upper semi-continuity in $L^2(\mathbb{R}^N)$ at $\varepsilon=0$. In this section, we strengthen this
study and  prove that
the upper semi-continuity of $\mathcal{A}_\varepsilon$  may happen in $H^1(\mathbb{R}^N)$ at $\varepsilon=0$.

For the upper semi-continuity, we  also give a  further assumption as in \cite{Wang1}, that is, $f$ satisfies that for all $x\in\mathbb{R}^N$ and
$s\in\mathbb{R}$,
\begin{align}\label{a6}
 \Big|\frac{\partial}{\partial s}f(x,s)\Big|\leq \alpha_4|s|^{p-2}+\psi_4(x),
\end{align}
where $\alpha_4>0$,  $\psi_4\in L^\infty(\mathbb{R}^N)$ if $p=2$ and $\psi_4\in L^{\frac{p}{p-2}}(\mathbb{R}^N)$ if $p>2$.

Let $\varphi_0$ be the continuous cocycle associated with the problem (\ref{eq1})-(\ref{eq2}) for $\varepsilon=0$. That is to say,
$\varphi_0$ is a deterministic non-autonomous cocycle over $\mathbb{R}$. Denote by $\mathcal{D}_0$ the collection of some families of deterministic nonempty subsets of $L^2(\mathbb{R}^N)$:
$$
\mathcal{D}_0=\{B=\{B(\tau)\subseteq L^2(\mathbb{R}^N);\tau\in \mathbb{R}\}; \lim\limits_{t\rightarrow+\infty}e^{-\delta t}\|B(\tau-t)\|=0,\tau\in\mathbb{R},\delta<\lambda\},
$$ where $\lambda$ is as in (\ref{pr1}).
As a special case of Theorem 4.9, under the  assumptions (\ref{a1})-(\ref{a5}), $\varphi_0$ has a common $\mathcal{D}_0$-pullback attractor $\mathcal{A}_0=\{ \mathcal{A}_0(\tau);\tau\in\mathbb{R}\}$ in both $ L^2(\mathbb{R}^N)$ and $H^1(\mathbb{R}^N)$.

To prove the upper semi-continuity of $\mathcal{A}_\varepsilon$ at $\varepsilon=0$, we have to check that the conditions (\ref{ff8})-(\ref{ff12}) in Theorem 2.8 hold in
$L^2(\mathbb{R}^N)$ and $H^1(\mathbb{R}^N)$ point by point. But (\ref{ff8})-(\ref{ff11}) have been achieved, see Corollary 7.2, Lemma 7.5 and equality (7.31) in \cite{Wang1}. We only need to prove the condition (\ref{ff12}) holds in
$H^1(\mathbb{R}^N)$.\\

\textbf{ Lemma  5.1.} \emph{ Assume that (\ref{a1})-(\ref{a5}) hold. Then for every  $\tau\in\mathbb{R}$ and $\omega\in \Omega$, the union $\cup_{\varepsilon\in(0,1]}\mathcal{A}_{\varepsilon}(\tau,\omega)$ is precompact in $H^1(\mathbb{R}^N)$.}\\

\emph{Proof}\  \ For any $\epsilon>0$, it suffices to show that for  every fixed $\tau\in\mathbb{R}$ and $\omega\in \Omega$, the set $\cup_{\varepsilon\in(0,1]}\mathcal{A}_{\varepsilon}(\tau,\omega)$
has finite $\epsilon$-nets in  $H^1(\mathbb{R}^N)$. Let $\chi=\chi(\tau,\omega)\in \cup_{\varepsilon\in(0,1]}\mathcal{A}_{\varepsilon}(\tau,\omega)$. Then there exists a $\varepsilon\in (0,1]$ such that
 $\chi(\tau,\omega)\in \mathcal{A}_{\varepsilon}(\tau,\omega)$. By the invariance of $\mathcal{A}_{\varepsilon}(\tau,\omega)$,
it follows that there is a $u_0\in \mathcal{A}_{\varepsilon}(\tau-t,\vartheta_{-t}\omega)$ such that
\begin{align}\label{mm00}
\chi(\tau,\omega)=\varphi_{\varepsilon}(t,\tau-t,\vartheta_{-t}\omega,u_0)=u_\varepsilon(\tau,\tau-t,\vartheta_{-\tau}\omega,u_0)\ \ (\mbox{by}\ (\ref{eq00})),
\end{align} for all $t\geq0$.
 Give $R>0$, denote by $\mathcal{O}^c_{R}=\mathbb{R}^N-\mathcal{O}_{R}$, where $\mathcal{O}_R=\{x\in \mathbb{R}^N; |x|\leq R\}$.
 Note that $\mathcal{A}_{\varepsilon}(\tau,\omega)\in\mathcal{D}$. Then by Lemma 4.4, for every $\epsilon>0$, there exist $T=T(\tau,\omega,\epsilon)\geq2$ and $R=R(\tau,\omega,\epsilon)>1$ such that the solution $u$ of problem (\ref{eq1})-(\ref{eq2})  satisfies
 for all $t\geq T$,
\begin{align}\label{mm01}
\|u_\varepsilon(\tau,\tau-t,\vartheta_{-\tau}\omega,u_0)\|_{H^1(\mathcal{O}^c_R)}\leq \epsilon.
\end{align}
Then by (\ref{mm00})-(\ref{mm01}), we have
\begin{align}\label{mm02}
\|\chi(\tau,\omega)\|_{H^1(\mathcal{O}^c_R)}\leq \epsilon,\  \ \ \mbox{for all}\ \chi\in\cup_{\varepsilon\in(0,1]}\mathcal{A}_{\varepsilon}(\tau,\omega).
\end{align}
On the other hand, by Lemma 4.6, there
exist a projector $P_{N_0}$ and  a $T=T(\tau,\omega,\epsilon)\geq2$ such that for all $t\geq T$
\begin{align} \label{mm03}
 \| (I-P_{N_0})\tilde{u}_\varepsilon(\tau,\tau-t,\vartheta_{-\tau}\omega, \tilde{u}_0)\|_{H^1_0(\mathcal{O}_{R\sqrt{2}})}\leq\epsilon,
\end{align}  where $\tilde{u}_\varepsilon$ is the cut-off of $u_\varepsilon$ on the domain $\mathcal{O}_{R\sqrt{2}}$, by (\ref{vuu}). Because
$P_{N_0}\tilde{u}_\varepsilon\in H_{N_0}$,
where $H_{N_0}=\mbox{span}\{e_1,_2,...,e_{N_0}\}$ is a finite dimension space and $P_{N_0}\tilde{u}_\varepsilon(\tau,\tau-t,\vartheta_{-\tau}\omega, \tilde{u}_0)$ is bounded
in $H_{N_0}$ which is compact.¡¡Therefore there exist some finite points  $v_1,v_2,...,v_s\in H_{N_0}$ such that
\begin{align} \label{mm04}
 \| P_{N_0}\tilde{u}_\varepsilon(\tau,\tau-t,\vartheta_{-\tau}\omega, \tilde{u}_0)-v_i\|_{H^1_0(\mathcal{O}_{R\sqrt{2}})}\leq\epsilon.
\end{align}
Thus by (\ref{mm00}), (\ref{mm03}) and (\ref{mm04}) are rewrote as
\begin{align} \label{mm05}
 \| (I-P_{N_0})\chi(\tau,\omega)\|_{H^1_0(\mathcal{O}_{R\sqrt{2}})}\leq\epsilon,\ \mbox{and}\ \| P_{N_0}\chi(\tau,\omega)-v_i\|_{H^1_0(\mathcal{O}_{R\sqrt{2}})}\leq\epsilon,
\end{align} for all $\chi\in\cup_{\varepsilon\in(0,1]}\mathcal{A}_{\varepsilon}(\tau,\omega)$. We now define $\tilde{v}_i=\tilde{v}_i(x)=0$ if $x\in \mathcal{O}^c_{R\sqrt{2}}$ and $\tilde{v}_i=v_i$ if $x\in \mathcal{O}_{R\sqrt{2}}$. Then for every $i=1,2,...,s$, $\tilde{v}_i\in H^1(\mathbb{R}^N)$. Furthermore, by (\ref{mm02}) and (\ref{mm05}), we have
\begin{align} \label{}
 \|\chi(\tau,\omega)-\tilde{v}_i\|_{H^1(\mathbb{R}^N)}&\leq \|\chi(\tau,\omega)-\tilde{v}_i\|_{H^1(\mathcal{O}^c_{R\sqrt{2}})}+\|\chi(\tau,\omega)-\tilde{v}_i\|_{H^1_0(\mathcal{O}_{R\sqrt{2}})}\notag\\
 &\leq\|\chi(\tau,\omega)\|_{H^1(\mathcal{O}^c_{R\sqrt{2}})}+\|P_N\chi(\tau,\omega)-\tilde{v}_i\|_{H^1_0(\mathcal{O}_{R\sqrt{2}})}\notag\\
 &
 +\|(I-P_{N_0})\chi(\tau,\omega)\|_{H^1_0(\mathcal{O}_{R\sqrt{2}})}\leq 3\epsilon,\notag
\end{align} for all $\chi\in\cup_{\varepsilon\in(0,1]}\mathcal{A}_{\varepsilon}(\tau,\omega)$. Thus
$\cup_{\varepsilon\in (0,1]}\mathcal{A}_\varepsilon(\tau,\omega)$ has finite $\epsilon$-nets in $H^1(\mathbb{R}^N)$, which implies that the union $\cup_{\varepsilon\in(0,1]}\mathcal{A}_{\varepsilon}(\tau,\omega)$ is precompact in $H^1(\mathbb{R}^N)$.
$\ \ \  \ \ \ \ \ \ \Box$\\

We  then obtain that the family of random attractors $\mathcal{A}_\varepsilon$ indexed by $\varepsilon$ converges to the deterministic $\mathcal{A}_0$ in $H^1(\mathbb{R}^N)$ in the following sense,\\

\textbf{Theorem 5.2.} \emph{ Assume that (\ref{a1})-(\ref{a5}) and (\ref{a6}) hold. Then for each $\tau\in\mathbb{R}$ and $\omega\in \Omega$,
$$
\lim\limits_{\varepsilon\downarrow0}\mbox{dist}_{H^1}(\mathcal{A}_\varepsilon(\tau,\omega), \mathcal{A}_0(\tau))=0
$$ where $\mbox{dist}_{H^1}$ is the Haustorff semi-metric in $H^1(\mathbb{R}^N)$.}\\

\section{Existence of random equilibria for the generated cocycle}

It is known that the random equilibrium is  a special case of omega-limit sets. The corresponding notion in deterministic case is fixed points or stationary solutions.
We can refer to \cite{Arn,Chues} for the definitions and applications.
The problem of the construction of equilibria for a general random dynamical system is rather complicate \cite{Chues}. Recently, \cite{Zhao1,Zhao2} obtained the existence of unique random equilibrium for stochastic reaction-diffusion equation with autonomous term on bounded domains or a unbounded Poincar\'{e} domains. Gu \cite{Gu1} proved that the stochastic FitzHugh¨CNagumo lattice equations driven by fractional Brownian motions possesses a unique equilibrium.

However, we here introduce the random equilibrium under the circumstance of non-autonomous stochastic dynamical system. In particular, we have\\

\textbf{Definition 6.1.} \emph{Let $(\Omega,\mathcal {F},{P},\{\vartheta_t\}_{t\in\mathbb{R}})$ be a measurable dynamical system. A random variable $u^*: \mathbb{R}\times\Omega\mapsto X$
is said to be an equilibrium (or fixed point, or stationary solution) of the cocycle $\varphi$ if it is invariant under $\varphi$, i.e., if
$$
\varphi(t,\tau,\omega, u^*(\tau,\omega))=u^*(\tau+t,\vartheta_t\omega)\ \  for\ all\ t\geq0,\ \tau\in\mathbb{R},\ \omega\in\Omega.
$$}

In this paper, we will prove the existence of equilibrium
for stochastic non-autonomous reaction-diffusion equation on the whole space $\mathbb{R}^N$.
We assume the coefficient $\lambda>\alpha_3$, where $\alpha_3$ is as in (\ref{a3}) and $\lambda$ is as in (\ref{pr1}). For convenience, here we write $\varepsilon=1$.
 First, we have\\

\textbf{Lemma 6.2.} \emph{Suppose that $g\in L^2(\mathbb{R}^N)$, $f$ and $g$ satisfies (\ref{a1})-(\ref{a5}) and $\lambda>\alpha_3$. Let the initial values $u_{0,i}=u(\tau-t_i,\vartheta_{-\tau}\omega)(i=1,2),t_1<t_2$. Then there exists a constant $b_0$ such that
the solution of problem (1.1)with initial value $u_{0,i}$ satisfies the following decay  property:
\begin{align} \label{}
& \ \  \ \ \ \ \ \ \ \|u(\tau, \tau-t_1,\vartheta_{-\tau}\omega, u_{0,1})-u(\tau, \tau-t_2,\vartheta_{-\tau}\omega, u_{0,2})\|^2\leq\notag\\
&2\Big(e^{-b_0t_1}e^{-2\omega(t_1)}\|u_{0,1}\|^2+2e^{-b_0 t_2}e^{-2\omega(t_2)}\|u_{0,2}\|^2\Big)
+ce^{(b_0-b)t_1}\int_{-\infty}^{0}e^{\lambda s}z^2(s,\omega)(\|g(s+\tau,.)\|^2+1)ds,\notag
\end{align}where $c$ is a deterministic non-random constant.}\\

\emph{Proof}\ \  Put $\bar{v}=v(\tau, \tau-t_1,\vartheta_{-\tau}\omega, v_{0,1})-v(\tau, \tau-t_2,\vartheta_{-\tau}\omega, v_{0,2})$. Then by (\ref{pr1})  we have
\begin{align} \label{equ1}
 \frac{d}{dt}\|\bar{v}\|^2+b\|\bar{v}\|^2\leq 0,
\end{align} where $b=\lambda-\alpha_3$. By applying Gronwall lemma to (6.1) over the interval $[\tau-t_1,\tau]$, we immediately get
\begin{align} \label{equ2}
 |\bar{v}(\tau)\|^2 &\leq e^{-bt_1}\|v(\tau-t_1, \tau-t_2,\vartheta_{-\tau}\omega, v_{0,2})-v_{0,1}\|^2\notag\\
 &\leq 2e^{-bt_1}\|v(\tau-t_1, \tau-t_2,\vartheta_{-\tau}\omega, v_{0,2})\|^2+2e^{-bt_1}\|v_{0,1}\|^2.
\end{align}
Choose
\begin{align} \label{equ3}
0<b_0<b.
\end{align}
By (\ref{se01}),  we have
\begin{align} \label{equ4}
\frac{d}{dt}\|v\|^2+b_0\|v\|^2\leq c z^2(t,\omega)(\|g(t,.)\|^2+1).
\end{align}
Then by Gronwall lemma again, we find that
\begin{align} \label{}
 &\|v(\tau-t_1, \tau-t_2,\vartheta_{-\tau}\omega, v_{0,2})\|^2\notag\\
 &\leq e^{b_0(t_1-t_2)}\|v_{0,2}\|^2+c\int_{\tau-t_2}^{\tau-t_1}e^{-b_0(\tau-t_1-s)}z^2(s,\vartheta_{-\tau}\omega)(\|g(s,.)\|^2+1)ds\notag\\
 &\leq  e^{b_0(t_1-t_2)}\|v_{0,2}\|^2+c e^{b_0 t_1}e^{2\omega(-\tau)}\int_{-\infty}^{0}e^{b_0 s}z^2(s,\omega)(\|g(s+\tau,.)\|^2+1)ds,\notag
\end{align}
from which and (\ref{equ2}) it follows that
\begin{align} \label{equ5}
 |\bar{v}(\tau)\|^2&\leq 2e^{-bt_1}\|v(\tau-t_1, \tau-t_2,\vartheta_{-\tau}\omega, v_{0,2})\|^2+2e^{-bt_1}\|v_{0,1}\|^2\notag\\
 &\leq2e^{-bt_1}\|v_{0,1}\|^2+2e^{(b_0-b)t_1}e^{-b_0 t_2}\|v_{0,2}\|^2\notag\\
 &+ce^{(b_0-b)t_1}e^{2\omega(-\tau)}\int_{-\infty}^{0}e^{\lambda s}z^2(s,\omega)(\|g(s+\tau,.)\|^2+1)ds\notag\\
 &\leq 2\Big(e^{-b_0t_1}\|v_{0,1}\|^2+e^{-b_0 t_2}\|v_{0,2}\|^2\Big)\notag\\
 &+ce^{(b_0-b)t_1}e^{2\omega(-\tau)}\int_{-\infty}^{0}e^{\lambda s}z^2(s,\omega)(\|g(s+\tau,.)\|^2+1)ds,
\end{align}
where we have used $e^{(b_0-b)t_1}\leq 1$ for $b_0<b$. By the equality $v(t)=z(t,\omega)u(t)=e^{-\omega(t)}u(t)$,  we get
\begin{align} \label{equ5}
 |\bar{u}(\tau)\|^2&\leq
 e^{-2\omega(-\tau)}|\bar{v}(\tau)\|^2\notag\\
 &\leq 2e^{-2\omega(-\tau)} \Big(e^{-b_0t_1}\|v_{0,1}\|^2+2e^{-b_0 t_2}\|v_{0,2}\|^2\Big)\notag\\
 &+ce^{(b_0-b)t_1}\int_{-\infty}^{0}e^{\lambda s}z^2(s,\omega)(\|g(s+\tau,.)\|^2+1)ds\notag\\
 &=2e^{-2\omega(-\tau)} \Big(e^{-b_0t_1}z^2(\tau-t_1,\vartheta_{-\tau}\omega)\|u_{0,1}\|^2+2e^{-b_0 t_2}z^2(\tau-t_2,\vartheta_{-\tau}\omega)\|u_{0,2}\|^2\Big)\notag\\
 &+ce^{(b_0-b)t_1}\int_{-\infty}^{0}e^{\lambda s}z^2(s,\omega)(\|g(s+\tau,.)\|^2+1)ds\notag\\
 &=2\Big(e^{-b_0t_1}e^{-2\omega(t_1)}\|u_{0,1}\|^2+2e^{-b_0 t_2}e^{-2\omega(t_2)}\|u_{0,2}\|^2\Big)\notag\\
 &+ce^{(b_0-b)t_1}\int_{-\infty}^{0}e^{\lambda s}z^2(s,\omega)(\|g(s+\tau,.)\|^2+1)ds,
\end{align} which finishes the proof.$\ \ \  \ \ \ \ \ \ \Box$\\

According to Lemma 6.2, we set $\lambda_0< b_0$ and  define the collection $\mathcal{D}$ by
$$
\mathcal{D}=\{B=\{B(\tau,\omega); \tau\in\mathbb{R},\omega\in\Omega\}; \lim\limits_{t\rightarrow+\infty}e^{-\lambda_0t}\|B(\tau-t,\vartheta_{-t})\omega\|^2=0,\ \mbox{for}\ \tau\in\mathbb{R},\omega\in\Omega\}.
$$ Then we have  the convergence result about the solution of problem (1.1)-(1.2) in $L^2(\mathbb{R}^N)$.\\

\textbf{Lemma 6.3.} \emph{Suppose that $g\in L^2(\mathbb{R}^N)$, $f$ and $g$ satisfies (\ref{a1})-(\ref{a5}) and $\lambda>\alpha_3$. Let $ B=\{B(\tau,\omega);\tau\in\mathbb{R},\omega\in\Omega\}\in\mathcal{D}$. Then for $\tau\in\mathbb{R},\omega\in\Omega$, there exists  a unique element $u^*=u^*(\tau,\omega)\in L^2(\mathbb{R}^N)$ such that
\begin{align} \label{}
\lim\limits_{t\rightarrow+\infty} u(\tau, \tau-t,\vartheta_{-\tau}\omega, u_{0})=u^*(\tau,\omega),\ \  \mbox{in}\ L^2(\mathbb{R}^N),\notag
\end{align} where $u_0\in B(\tau-t,\vartheta_{-t}\omega)$.
Furthermore, the convergence is uniform (w.r.t $u_0\in B(\tau-t,\vartheta_{-t}\omega)$).}\\

\emph{Proof}\ \ If $u_{0,i}\in B(\tau-t_i,\vartheta_{-t_i}\omega)$, then we have
 $\lim\limits_{t_i\rightarrow+\infty}e^{-b_0t_i}e^{-2\omega(t_i)}\|u_{0,i}\|^2=0$ for $i=1,2$. Thus the result is derived directly from Lemma 6.2.
 $\ \ \  \ \ \ \ \ \ \Box$\\

\textbf{Lemma 6.4.} \emph{Suppose that $g\in L^2(\mathbb{R}^N)$, $f$ and $g$ satisfies (\ref{a1})-(\ref{a5}) and $\lambda>\alpha_3$. Then for $\tau\in\mathbb{R},\omega\in\Omega$, the  element $u^*=u^*(\tau,\omega)$ defined in Lemma 6.3 is a unique random equilibria for the cocycle $\varphi$  defined by (\ref{eq0}) in $L^2(\mathbb{R}^N)$, i.e.,
\begin{align} \label{}
\varphi(t, \tau,\omega, u^*(\tau,\omega))=u^*(\tau+t,\vartheta_{t}\omega),\ \ \mbox{for every }\ t\geq0,\ \tau\in\mathbb{R},\ \omega\in\Omega.\notag
\end{align}
Furthermore, the random equilibria $\{u^*(\tau,\omega), \tau\in\mathbb{R},\ \omega\in\Omega\}$  is the unique element of the pullback attractor $\mathcal{A}=\{\mathcal{A}(\tau,\omega); \tau\in\mathbb{R},\ \omega\in\Omega\}$ for the cocycle $\varphi$, i.e., for every $\tau\in\mathbb{R},\omega\in\Omega$,
$\mathcal{A}(\tau,\omega)=\{u^*(\tau,\omega)\}.$}\\

\emph{Proof}\ \  By the definition of the cocycle, $\varphi(t,\tau-t, \vartheta_{-t}\omega, u_0)=u(\tau,\tau-t,\vartheta_{-\tau}\omega,u_0)$, then for for every $\tau\in\mathbb{R},\omega\in\Omega$, we have
\begin{align} \label{6.8}
u^*(\tau,\omega)=\lim\limits_{t\rightarrow+\infty} \varphi(t, \tau-t,\vartheta_{-t}\omega, u_{0}),
\end{align} where $u_0\in B(\tau-t,\vartheta_{-t}\omega)$. Thus by the continuity and the cocycle property  of $\varphi$ and (\ref{6.8}), we find that for every $t\geq0,\tau\in\mathbb{R},\omega\in\Omega$,
\begin{align} \label{}
\varphi(t, \tau,\omega, u^*(\tau,\omega))&=\varphi(t, \tau,\omega, .)\circ \lim\limits_{s\rightarrow+\infty} \varphi(s, \tau-s,\vartheta_{-s}\omega, u_{0})\notag\\
&=\lim\limits_{t\rightarrow+\infty} \varphi(t, \tau,\omega, .)\circ \varphi(s, \tau-s,\vartheta_{-s}\omega, u_{0})\notag\\
&=\lim\limits_{t\rightarrow+\infty} \varphi(t+s, \tau-s,\vartheta_{-s}\omega, u_{0})\notag\\
&=\lim\limits_{t\rightarrow+\infty} \varphi(t+s, (\tau+t)-t-s, \vartheta_{-s-t}\vartheta_{t}\omega, u_0)\notag\\
&=u^*(\tau+t,\vartheta_{t}\omega),\notag
\end{align}
which also implies  the invariance of $\mathcal{A}$, that is, $\varphi(t, \tau,\omega, \mathcal{A}(\tau,\omega))=\mathcal{A}(\tau+t,\vartheta_t\omega)$. The compactness of
$\mathcal{A}(\tau,\omega)$ is obvious and the attracting property follows from (\ref{6.8}). $\ \ \  \ \ \ \ \ \ \Box$\\

\textbf{Remark 6.5.}\ \  We notice that  by Theorem 4.9, the equilibria $u^*\in H^1(\mathbb{R}^N)$. In particular, we further have $u^*\in L^p(\mathbb{R}^N)$.\\

\textbf{Acknowledgments:}

 This work was supported by Chongqing Basis and Frontier Research Project of China (no. cstc2014jcyjA00035) and National Natural Science Foundation of China
(no.11071199).\\

\end{document}